\documentclass[12pt]{article}
\usepackage{amsfonts,amsmath,color,graphicx,subfigure,bm,paralist,amsthm}

\usepackage{mathrsfs}

\usepackage{xcolor}

\usepackage{tabularx}
\title{Automatic Decoupling and  Index-aware  \colp{Model-Order Reduction} for Nonlinear \colp{Differential-Algebraic Equations}}

\author{N. Banagaaya \thanks{Max Planck Institute for 
Dynamics of Complex Technical Systems, Magdeburg, Germany
        ({\tt \{banagaaya, grundel, benner\}@mpi-magdeburg.mpg.de}).}
        \and
        G. Al\`{\i}\thanks{Dept. of Mathematics, University of Calabria,
        and INF-Group c. Costanza, 
        Cataract DI Rene (Costanza), Italy
       ({\tt Giuseppe.ali@unical.it}).}\and  S. Grundel\footnotemark[1] \and P. Benner \footnotemark[1]
        }
        
\usepackage{geometry}
\def\E{{{\bf E}}}
\def\A{{{\bf A }}}
\def\B{{{\bf B}}}
\def\C{{{\bf C}}}

\def\U{{{\bf U}}}

\def\d{{{\bf d}}}

\def\g{{{\bf g}}}

\def\M{{{\bf M}}}
\def\N{{{\bf N}}}

\def\J{{{\bf J}}}
\def\W{{{\bf W}}}

\def\RR{{{\mathbb R}}}

\def\Nul{{{\mathrm {Ker} }}}
\def\Im{{{\mathrm {Im} }}}

\def\f{{{\bf f}}}
\def\x{{{\bm x}}}
\def\u{{{\bm u}}}
\def\y{{{\bm y}}}

\def\Q{{{\bf Q}}}
\def\P{{{\bf P}}}

\def\q{{{\bf q}}}
\def\p{{{\bf p}}}

\def\s{{{\bf s}}}

\def\V{{{\bf V}}}

\newcommand{\tr}{\mathrm{T}}

\theoremstyle{definition}
\newtheorem{example}{Example}[section]
\newtheorem{definition}{Definition}[section]

\newcommand{\colb}[1]{{ #1}}
\newcommand{\colr}[1]{{#1}}
\newcommand{\colg}[1]{{ #1}}
\newcommand{\colp}[1]{{ #1}}
\date{}
\begin{document}
\maketitle
\abstract{We  extend   the index-aware \colp{model-order reduction}    method to \colp{systems of} nonlinear differential-algebraic equations   with a special nonlinear term $\f(\E\x),$ where $\E$ is a singular matrix.
Such  nonlinear \colp{differential-algebraic equations}   arise\colb{,} for example\colp{,} in the  spatial discretization of the gas flow in pipeline  networks.
In practice, mathematical models of real-life processes pose challenges when used in numerical simulations, due to complexity and system size.
Model-order reduction  aims to eliminate this problem   by generating reduced-order
models  that have lower   computational cost  to simulate, yet accurately represent
the original large-scale system behavior. However, direct reduction and simulation  of  nonlinear \colp{differential-algebraic equations}  is difficult due to hidden constraints which affect the choice of numerical integration  methods  and \colp{model-order reduction} techniques.
We propose an extension of \colp{index-aware model-order reduction} methods to nonlinear \colp{differential-algebraic equations}  without any kind of linearization.
The proposed \colp{model-order reduction} approach involves
automatic decoupling of nonlinear \colp{differential-algebraic equations}  into nonlinear \colp{ordinary differential equations} and algebraic
equations.  This allows applying  standard \colp{model-order reduction} techniques 
\colp{to} both parts without worrying about the index. The same procedure can also be used to simulate nonlinear \colp{differential-algebraic equations using standard integration schemes}. 
We illustrate the performance of our proposed method \colp{for} nonlinear \colp{differential-algebraic equations}  arising   from   gas flow models in  pipeline   networks.}
\section{Introduction}
We consider nonlinear  \colp{differential-algebraic equations (}DAEs\colp{)} of the form;
\begin{subequations}
\label{Eqn:1}
\begin{align}
\label{Eqn:1a}
 \E\x'&=\A\x+\f(\E\x)+\B\u,\quad\mathbf{E}\mathbf{x}(0)=\mathbf{E}\mathbf{x}_0, \\ 
 \y&= \C\x,\label{Eqn:1b}
\end{align}
\end{subequations}
where $\f(\E\x)\in \mathbb{R}^{n}$ and  $\E$ is a singular matrix,  $\mathbf{A}\in  \mathbb{R}^{n\times n},\, \mathbf{B} \in  \mathbb{R}^{n\times m},\,  \mathbf{C}\in  \mathbb{R}^{\ell \times n}$. \colb{The symbol $'$ denote  time differentiation.}
  $\x\in \mathbb{R}^{n}$ and   $\mathbf{y}\in \mathbb{R}^{\ell}$  are the state  and output vectors, respectively.  The input function   ${\u}\in \mathbb{R}^{m}$ must be  smooth enough\colp{, with the smoothness requirements} depending on  the index of the DAE.
DAEs are known to be difficult to simulate and the  level of difficulty is  measured using index concepts such as differential index, tractability index, etc. The higher the index, the more difficult to simulate the DAE.
Moreover, in practice, \colp{often} such \colp{\emph{descriptor} } systems have very large  $n,$
compared to  the number $m$ of inputs and the number $\ell$  of outputs, which  are typically small.  Despite the  ever increasing computational power, dynamic  simulation using  the  system \eqref{Eqn:1} is costly, see \cite{morGruHR16,morAnt05,morBauBF14}.
We are interested in a fast and stable prediction of the dynamics of DAE models,
and therefore the application of  \colp{model-order reduction (}MOR\colp{)} is vital.
MOR aims
to reduce the computational burden by generating  \colp{reduced-order models (}ROMs\colp{)} that have   lower computational cost
to simulate, yet accurately represent the original  large-scale system  behavior. MOR replaces \eqref{Eqn:1} by a ROM
\begin{subequations}
 \label{Eqn:Red}
 \begin{align}
\mathbf{E}_r\mathbf{x}_r'&=\mathbf{A}_r\mathbf{x}_r+\mathbf{f}_r(\E_r\mathbf{x}_r)+\mathbf{B}_r\mathbf{u},\quad\mathbf{E}_r\mathbf{x}_r(0)=\mathbf{E}_r\mathbf{x}_{r_0},\\
\mathbf{y}_r&=\mathbf{C}_r\mathbf{x}_r,
\end{align}
\end{subequations}
where $\mathbf{E}_r,\mathbf{A}_r\in \mathbb{R}^{r\times r},\, \mathbf{f}_r(\E_r\mathbf{x}_r)\in \mathbb{R}^{r},\mathbf{B}_r\in\mathbb{R}^{r\times m}$ and $\mathbf{y}_r\in\mathbb{R}^{\ell},
\mathbf{C}_r\in \mathbb{R}^{\ell\times r},$ such that the reduced-order of the state vector $\mathbf{x}_r\in\mathbb{R}^{r}$ is $r\ll n.$
A good ROM should have small approximation error
$\Vert  \mathbf{y}-\mathbf{y}_r\Vert $ in a suitable norm $\Vert  .\Vert $ for a desired range of inputs $\mathbf{u}.$
There exist many
MOR methods for nonlinear  systems such as \colp{proper orthogonal decomposition (}POD\colp{)},  \colp{  POD in conjunction with the discrete empirical interpolation method} (POD-DEIM), see \cite{morBauBF14}.  However, applying these MOR methods directly to DAEs  leads to ROMs which are  inaccurate or  very difficult to simulate and sometimes have no solution, see \cite{morGruHR16,BanaBook}.
It is a common practice to first convert  nonlinear  DAEs to \colp{ordinary differential equations (}ODEs\colp{)} by using index reduction (reformulation) techniques in order to be  able to apply standard MOR methods for nonlinear systems such as POD.
 However, this  index reduction may lead  to drift-off effects or instabilities in the numerical solutions and may also  depend on the structure of the nonlinear DAE.
 In \cite{Bana:2014}, IMOR methods were proposed  to eliminate the index problem to allow employing standard techniques with ease.   However, \colp{these} method\colp{s} w\colp{ere} dedicated to linear DAEs.
We propose an  index-aware MOR (IMOR) \colp{method} for nonlinear DAEs \colb{of the form \eqref{Eqn:1}} which does not involve any kind of linearization. This approach  is \colp{realized} in two steps. The  first step involves  automatically decoupling  the nonlinear DAEs  into nonlinear differential and algebraic parts.
Then, each part can be reduced separately using standard MOR techniques. 
The decoupled  system generated from the first step can be used for  numerical simulations by applying numerical integration on the ODE part and then solving the algebraic part. 


The paper is organized as follows. In {S}ection \ref{sec:1}, we discuss the background of  decoupling of DAEs and the tractability index.   In {S}ection \ref{sec:2}, we propose   the automatic  decoupling of nonlinear DAEs of the form \eqref{Eqn:1} using special projectors.
In {S}ection \ref{sec:4}, we discuss the proposed 
IMOR method
for nonlinear DAEs.   In {S}ection \ref{sec:5}, we  apply the proposed IMOR method to the nonlinear DAEs arising from  gas transport networks.
In the final section, we present some numerical examples
illustrating the {performance} of the proposed method.
\section{Decoupling of \colg{linear constant coefficient} DAEs}
\label{sec:1}
In this section, we repeat the procedure of decoupling linear DAEs and the theory it is based on, as this is the basis for the nonlinear decoupling.
\subsection{ Weierstra{\ss}  canonical form}
\colp{Our}  \colp{d}ecoupling strategy   was initially used to 
understand the   underlying structure of linear constant coefficient DAEs via \colp{the}  Weierstra{\ss} canonical form \colg{\cite{Canonical}}.   Assuming $\colb{\f(\E\x)}=0$  and \colp{that the} matrix pencil $(\E,\A)$ is regular\colp{,} \eqref{Eqn:1} can be written as a  Weierstra\ss-Kronecker canonical form which leads to an 
equivalent decoupled system
\begin{subequations}
\label{Eqn:krock}
\begin{align}
  \tilde{\x}'_1 &= \J \tilde{\x}_1  + \tilde{\B}_1\u, \quad \tilde{\x}_1(0)= \tilde{\x}_{1_0}, \label{Eqn:krocka}\\
  \tilde{\x}_2& = - \sum_{i=0}^{\mu-1}\N^i \tilde{\B}_2\u^{(i)}, \label{Eqn:krockb}
\end{align}
\end{subequations}
where $\J \in \RR^{k \times k}$  and $\N \in \RR^{(n-k)\times (n-k) }$ is a nilpotent matrix with index $\mu.$ 
\colb{The vector} $\displaystyle \u^{(i)}=\frac{\mathrm{d}^i}{\mathrm{d}t^{\colb{i}}}\u\in  \RR^m$ is the $i$-th derivative 
of the input data. \colb{The control input matrices are }$\tilde{\B}_1\in \RR^{k\times m}$ \colb{and} $\tilde{\B}_2\in \RR^{(n-k)\times (n-k)}.$  Subsystems \eqref{Eqn:krocka} and \eqref{Eqn:krockb} represent the inherited ODE and algebraic part, 
 respectively\colb{,} and the solutions of \eqref{Eqn:1} can be obtained using $\x=\W\tilde{\x}$  where  $\W^{n\times n}$ is a nonsingular matrix  and $\displaystyle \tilde{\x}=\left(
                                                                                                                                                     \tilde{\x}_1^\tr,
                                                                                                                                                     \tilde{\x}_2^\tr
                                                                                                                                                    \right)^\tr.
 $ The vectors $ \tilde{\x}_1 \in\RR^{k}$ and $ \tilde{\x}_\colb{2}^{n-k}$ are commonly known as the slow and fast \colp{parts of the} solution, respectively.
An index concept  was introduced to classify different types of DAEs with respect to the difficulty arising in the theoretical and numerical treatment of a given DAE. 
\colp{``}Index\colp{''} is a notion used in the theory of DAEs for measuring the distance from a DAE to its related ODE. There are several definitions of a DAE index.
The index $\mu$ in \eqref{Eqn:krock}  is  known as  the Kronecker  \colp{(nilpotency, differentiability)} index. 
An equivalent decoupled system \eqref{Eqn:krock} shows the dependence of the solution of \colg{a linear} DAE   on the 
derivatives of the input function.  In  \eqref{Eqn:krockb}, we can observe that the input function has to be at least $\mu-1$ times differentiable.    The higher the index  the more differentiations of the input data are involved. 
 Since numerical differentiation is an unstable process, the index $\mu$ is a
measure of numerical difficulty when solving the DAE.
We can also observe that the initial condition $\tilde{\x}_1(0)$ of the differential part can be chosen arbitrary while  $\tilde{\x}_2(0)$ has to satisfy hidden constraint  $$\displaystyle \tilde{\x}_2(0) = - \sum_{i=0}^{\mu-1}\N^i \tilde{\B}_2\u^{(i)}(0).$$
Thus, DAE \eqref{Eqn:1} has a unique classical solutions if $\x(0)=\x_0$ is consistent. 
The index problem affects the choice of numerical integration schemes strongly if standard numerical integration schemes are applied \colp{to} DAEs directly without decoupling. This lead to the 
development of numerical integration schemes which were specifically \colb{designed} for DAEs, see \cite{Marz92,MehrBook}.  Hence, \colp{a promising} way to solve and apply MOR \colp{to} DAEs is to first split them into differential and algebraic parts, see \cite{BanaBook}.
According to \cite{Marz92}, transforming linear DAEs into a Kronecker canonical form is numerically infeasible and it is restricted to linear DAEs.
Due to this drawback other index concepts, such as  the tractability index, differentiation index, etc, see \cite{4lectures}, were proposed with  each of them stressing different aspects of the DAE.
\subsection{Tractability index}
In this paper, we consider the tractability index introduced in \cite{MarzBook}  and its generalization in \cite{Marz02} defined as in Definition \ref{defn1}.
\begin{definition}[Tractability index (\cite{Marz02}) ]
\label{defn1}
 Given a regular  matrix pair $(\E,\A).$  We define a matrix and projector chain by setting   ${\bf E}_{0} :={\bf E}$ and ${\bf A}_{0}:={\bf A}$, given by
 \begin{align}
 {\bf E}_{j+1}& := {\bf E}_{j}-{\bf A}_{j}{\bf Q}_{j},\quad {\bf A}_{j+1}:={\bf A}_{j}{\bf P}_{j},\quad \mbox{for}\, j\geq0,\label{proj}
\end{align}
 where ${\bf Q}_{j}$  are projectors onto  $\Nul\, {\bf E}_{j}$ and ${\bf P}_{j}=\mathrm{I}-{\bf Q}_{j}$.
 There exists an index $\gamma$ such that ${\bf E}_{\gamma}$ is non-singular and all ${\bf E}_{j}$ are singular for all $0\leq j < \gamma-1$. $\gamma$ is the \emph{tractability index} \colp{of a DAE}.
\end{definition}
This index criterion does not depend on the special choice of the projector functions $\Q_j$, see \cite{Marzlead}.
The tractability index has gained a lot of attention since  \colb{it}  can be calculated without the use of derivative arrays \cite{Marz_arrays2}.
Hence, it is numerical\colp{ly} feasible \colp{to compute the tractability index } compared  to  \colp{computing the} Kronecker index.
This is  the main tool in the decoupling of linear DAEs into their differential and algebraic parts, since it allows
  automatic decoupling procedure\colp{s}, see \cite{Marz02, Bana:2014,4lectures}. In order to decouple linear  DAEs with   index higher than one\colb{,}  so\colb{-}called canonical projectors were introduced  in \cite{Marz4} with  
  additional constraint\colb{s} $\mathbf{Q}_j\mathbf{Q}_i=0$, for $j\textgreater i.$ 
  Based on these projectors, special projector bases were introduced leading to a decoupled system of the same dimension as the original.
  A key step \colp{in} forming the projectors in \eqref{proj}
is to find the initial projectors $\mathbf{Q}_j$  spanning the
nullspaces of the usually sparse $\mathbf{E}_j$. Standard
ways of identifying the nullspace include  the singular value decomposition (SVD) or alike,
which do\colp{es} not utilize matrix patterns and can be expensive for
large-size matrices. \colp{One feasible way is} to  employ the sparse LU
decomposition-based routine, called LUQ, see \cite{Zhang}. This  routine was also used to construct the projector bases efficiently. This motivated the development of the index-aware MOR methods in \cite{Bana:2014}.
In \cite{Bana:2014},  equivalent  explicit and implicit decoupling \colp{methods}  were  proposed which are discussed \colp{in the following}. 
According to  \cite{Bana:2014}, if we also assume $\f(\colb{\E}\x)=0$ and \colp{that} the matrix pencil $(\E,\A)$ is regular, using \colp{the} special projectors proposed in \cite{Marz4} and projector bases introduced in \cite{Bana2014}\colp{,}
DAE  \eqref{Eqn:1} can be rewritten into an   equivalent explicit  decoupled system given by
 \begin{subequations}
\label{Eqn:10}
\begin{align}
\xi_{p}'&=\mathbf{A}_{p}\xi_{p}+\mathbf{B}_{p}\mathbf{u},\quad\xi_p(0)=\xi_{p_0},\label{Eqn:10a}\\
 \xi_{q} &=\sum_{j=0}^{\gamma-1}\mathcal{L}^{j}\left( \mathbf{A}_{q}\xi_{p}^{(j)}+\mathbf{B}_{q}\u^{(j)} \right)\label{Eqn:10b} \\
 \mathbf{y}&=\mathbf{C}_{p}\xi_{p}+\mathbf{C}_{q}\xi_{q},
\end{align}
\end{subequations}
{where}  $\mathcal{L}\in \mathbb{R}^{n_q\times n_q}$ is a nilpotent matrix with    index $\gamma.$
$\u^{(j)}$ and $\xi_{p}^{(j)}$ are the $j$-th derivatives with respect to $t$. 
  The subsystems \eqref{Eqn:10a} and  \eqref{Eqn:10b} correspond to the differential and algebraic parts  of system \eqref{Eqn:1}.  $\xi_{p}\in \mathbb{R}^{n_p}$ and $\xi_{q} \in \mathbb{R}^{n_q}$ are
 the differential and algebraic variables. The dimension of the decoupled system is given by $n=n_p+n_q.$  We can observe that decoupled system\colp{s}  \eqref{Eqn:krock} and  \eqref{Eqn:10} are equivalent if $\mathbf{A}_{q}=0.$
 Decoupled system \eqref{Eqn:10}  can be constructed
 automatically\colp{,} and \colp{thus} is numerically feasible. However, according to \cite{Bana:2014}  the decoupling procedure of \eqref{Eqn:10} involves  the inversion of non-singular  matrix ${\bf E}_{\gamma}$  which is costly for large-scale systems. 
 Then, the implicit version of \eqref{Eqn:10} was also proposed in \cite{Bana:2014} which  do\colp{es} not involve  the  inversion of  non-singular  matrix ${\bf E}_{\gamma}$.  Using this decoupling procedure\colp{,}  DAE  \eqref{Eqn:1} can be \colp{re-}written into an   equivalent implicit  decoupled system given by
 \begin{subequations}
\label{Eqn:11}
\begin{align}
\mathbf{E}_p\xi_{p}'&=\mathbf{A}_{p}\xi_{p}+\mathbf{B}_{p}\mathbf{u},\quad\xi_p(0)=\xi_{p_0},\label{Eqn:11a}\\
    \mathcal{L}_q\xi_{q} &=\sum_{j=0}^{\gamma-1}\mathbf{N}_q^{j}\left( \mathbf{A}_{q}\xi_{p}^{(j)}+\mathbf{B}_{q}\u^{(j)} \right) ,\label{Eqn:11b}\\
 \mathbf{y}&=\mathbf{C}_{p}\xi_{p}+\mathbf{C}_{q}\xi_{q},
\end{align}
\end{subequations}
 where $\mathbf{N}_q=\mathcal{L}\mathcal{L}_q^{-1}$ is also a nilpotent matrix with the same index $\gamma$ as $\mathcal{L}.$
The matrices   $\mathcal{L}_q\in \mathbb{R}^{n_q\times n_q}$ and $\mathbf{E}_p\in \mathbb{R}^{n_p\times n_p}$ are always    non-singular\colb{, see \cite{Bana2014}}.
  The subsystems \eqref{Eqn:11a} and  \eqref{Eqn:11b} correspond to the differential and algebraic parts  of system \eqref{Eqn:1}.  $\xi_{p}\in \mathbb{R}^{n_p}$ and $\xi_{q} \in \mathbb{R}^{n_q}$ are
 the differential and algebraic variables. \colb{We can observe that the inherited ODEs \eqref{Eqn:10a} and \eqref{Eqn:11a} of the explicit and implicit decoupled systems can be  simulated using standard ODE integration schemes. After obtaining the solutions
 of \eqref{Eqn:11a}, the
 algebraic part \eqref{Eqn:11b} can \colp{be} solved using numerical solvers such as    LU decomposition-based  routines.  }
 \colg{It is not straight forward to extend this to nonlinear DAEs. However, specific classes of nonlinear DAEs have been studied, }
usually those
appearing in practice, see \cite{Bana2018}.


\section{Decoupling of nonlinear DAEs  }
\label{sec:2}
In this section, we propose the  decoupling of a class of nonlinear DAEs of the form \eqref{Eqn:1a}. This decoupling strategy is an exten\colp{s}ion of  the  decoupling \colp{strategy for} linear DAEs proposed in \cite{Marz4}. 
\subsection{\colr{Decoupling using projectors}}
Assum\colp{e that} the tractability index  of \eqref{Eqn:1a} is independent of the nonlinearity, i.e., all projectors constructed using \colb{D}efinition \ref{defn1} are constant matrices.  Setting $\E_0=\E,\quad \A_0=\A,$   \eqref{Eqn:1a} can be written as 
\begin{align}
 \E_0\x'&=\A_0\x+\f(\E_0\x)+\B\u.\label{Eqn:2bb}
\end{align}
\colb{W}e choose \colb{a} projector $\Q_0$ such that $\Im \Q_0=\Nul \E_0 $ and its complementary projector $\P_0=\mathrm{I}-\Q_0.$ Using  \eqref{proj},  $$\E_1=\E_0-\A_0\Q_0, \quad  \A_1= \A_0\P_0,$$ which satisfy  the identities\colb{:}
\begin{align}
 \E_1\P_0&=\E_0,\quad \A_1-\E_1\Q_0=\A_0. 
\end{align}
Substituting the above identities into \eqref{Eqn:2bb} and simplifying leads to
\begin{align}
 \E_1\big[ \P_0\x'+\Q_0\x\big]&= \A_1\x+\f(\E_1\P_0\x)+\B\u. \label{Eqn:5b}
\end{align}
If we assume $\E_1$ \colp{to be} nonsingular, then \eqref{Eqn:5b} can be written as 
\begin{align}
 \P_0\x'+\Q_0\x&=  \E_1^{-1}\big[\A_1\x+\f(\E_1\P_0\x)+\B\u\big]. \label{Eqn:6}
\end{align}
Since $\E_1$ is nonsingular, then we say that   the nonlinear DAE \eqref{Eqn:1} is of tractability index $1.$ 
Left multiplying \eqref{Eqn:6} by  projectors $\P_0$ and $\Q_0$ separately,  we obtain the differential and algebraic subsystems, respectively\colp{,}  of  \eqref{Eqn:1} given by 
\begin{subequations}
\label{Eqn:decouple_index1}
\begin{align}
 \x_P'&=\P_0 \E_1^{-1}\A_0\x_P+ \P_0 \E_1^{-1}\f(\E_1\x_P) +\P_0 \E_1^{-1}\B\u,\quad \x_P(0)=\P_0\x(0),\\
 \x_Q&=\Q_0\E_1^{-1}\A_0\x_P+ \Q_0 \E_1^{-1}\f(\E_1\x_P) +\Q_0 \E_1^{-1}\B\u,\\
 \y&= \C\x_P+\C\x_Q,
\end{align}
\end{subequations}
where $\x_P=\P_0\x$ and $\x_Q=\Q_0\x.$
We can see that decoupled system \eqref{Eqn:decouple_index1} is of dimension $2n$ while the DAE \eqref{Eqn:1} is of dimension $n$. This implies that decoupling using projectors does not
preserve the dimension of the original DAE. In the next section, we discuss how to derive a decoupled system which preserves the dimension of the  nonlinear DAE \eqref{Eqn:1}.

\subsection{Explicit decoupling \colg{using bases} }
\label{Sec:Explicit}
Projector bases can be applied \colp{to} \eqref{Eqn:decouple_index1} as follows.
Let $n_q=\mathrm{dim}(\mathrm{Ker}\,{\bf E}_0)$ and  $n_q=n-n_p.$
If, we also  let  $\q_0 \in \Im \, \Q_0 $ and $\p_0\in \Im \, \P_0\colb{,}$  
\colb{t}hen, we can  expand $\bm{x}$ with respect to the bases, obtaining
\begin{align}
 \bm{x}=\q_0\xi_q+\p_0\xi_p,\label{soln:1}
\end{align}
where $\xi_q\in\mathbb{R}^{n_q}, \quad  \xi_p\in\mathbb{R}^{n_p},$
which implies that
$
 \bm{x}_P=\p_0\xi_p \quad \mbox{and} \quad  \bm{x}_Q=\q_0\xi_q$  in \eqref{Eqn:decouple_index1}.
\colb{T}he left inverses of column matrices  $\q_0\in \RR^{n\times n_q}$ and $\p_0\in \RR^{n\times n_p}$ are denoted by  $\q_0^{*\tr} \in \mathbb{R}^{n_q \times n}$ and $\p_0^{*\tr} \in \mathbb{R}^{n_p \times n},$  respectively.
 Substituting   $\bm{x}_P=\p_0\xi_p \, \mbox{and}$ \\ $\bm{x}_Q=\q_0\xi_q$  into \eqref{Eqn:decouple_index1}   leads to a  decoupled system \colg{which can be}
 \colg{l}eft multipl\colg{ied}   by the left inverses  $\p_0^{*\tr}$ and $\q_0^{*\tr},$ respectively\colp{.} \colp{This yields a } decoupled system  in  compact form:
\begin{subequations}
\label{eqn:303}
\begin{align}
\xi_p' &= {\bf A}_{p}\xi_p+\f_p(\xi_p)+{\bf B}_{p}\bm{u}, \quad \xi_p(0)=\p_0^{*\tr}\bm{x}(0), \label{eqn:303a}
\\
\xi_q &= {\bf A}_{q} \xi_p+\f_q(\xi_p)+{\bf B}_{q} \bm{u}, \label{eqn:303b}\\
 \y&= \C_p\xi_p+\C_q\xi_q,  \label{eqn:303c}
\end{align}
\end{subequations}
where 
\begin{align*}
 {\bf A}_{p}&=\p_0^{*\tr}{\bf E}_1^{-1}{\bf A}_{0}\p_0\in \RR^{n_p \times n_p},\, {\bf B}_{p}=\p_0^{*\tr}{\bf E}_1^{-1}{\bf B}\in \RR^{n_p \times m},\,  \C_p=\C\p_0 \in \RR^{\ell\times  n_p},\\
 \C_q&=\C\q_0 \in \RR^{\ell\times  n_q},\, {\bf A}_{q}=\q_0^{*\tr}{\bf E}_1^{-1}{\bf A}_{0}\p_0\in \RR^{n_q \times n_p} ,\, {\bf B}_{q}=\q_0^{*\tr}{\bf E}_1^{-1}{\bf B}\in \RR^{n_q \times m}.
\end{align*}
\colg{and} $$\f_p(\xi_p)=\p_0^{*\tr} \E_1^{-1}\f(\E_1\p_0\xi_p) \in \RR^{n_p}, \quad   \f_q(\xi_p)=\q_0^{*\tr} \E_1^{-1}\f(\E_1\p_0\xi_p)\in \RR^{n_q}.$$
We  can now  observe that the total dimension of the decoupled system is $n=n_p+n_q$, which is equal to the dimension of the nonlinear DAE \eqref{Eqn:1}. Instead of solving  the coupled nonlinear DAE \eqref{Eqn:1} we can now solve the 
decoupled nonlinear system \eqref{eqn:303}. We obtain the solution   $\xi_p$  by applying standard integration schemes to \eqref{eqn:303a}  and  the solutions of $\xi_q$ can be computed by post-processing using \eqref{eqn:303b}.
Then, the desire\colp{d} output solution can be obtained using \eqref{eqn:303c}. However, we can observe that  the coefficients of \eqref{eqn:303} involve computing the inverse of  $\E_1$  which is 
computationally expensive and requires  large storage  for large scale systems. Moreover, it also leads to   dense matrix coefficients of the decoupled system \eqref{eqn:303}.

\subsection{Implicit decoupling }
\label{sub:sec1}
In this subsection, we discuss \colp{a} decoupling strategy  which does not involve inversion of matrix $\E_1$. This is done as follows.
Substituting \eqref{soln:1} into \eqref{Eqn:5b}  leads to 
\begin{align}
 \begin{pmatrix}
   \E_1\p_0  & 0 
 \end{pmatrix}\begin{pmatrix}
 \xi_p\\
 \xi_q
 \end{pmatrix}'&=\begin{pmatrix}
   \A_0\p_0  & -\E_1\q_0 
 \end{pmatrix}\begin{pmatrix}
 \xi_p\\
 \xi_q
 \end{pmatrix}+\f(\E_1\p_0\xi_p)+\B\u. \label{Eqn:12}
\end{align}
Instead of inverting  matrix $\E_1,$ we can  decouple \eqref{Eqn:12} into differential and algebraic parts using column matrices   $\hat{\p}_0\in \mathbb{R}^{n\times n_p}$ and  $\hat{\q}_0 \in \mathbb{R}^{n\times n_p}$
proposed in \cite{Bana:2014} which are  defined as via   $\hat{\p}_0\in  \mathrm{Ker} \, \q_0^{\tr}\E_1^{\tr}$ and $\hat{\q}_0\in \mathrm{Ker} \, \p_0^{\tr}\E_1^{\tr}$.
\colp{L}eft multiplying \eqref{Eqn:12} by $\big( \hat{\p}_0^{\tr}\quad  \hat{\q}_0^{\tr} \big)^{\tr}$ lead\colp{s} to 
\begin{multline}
 \begin{pmatrix}
   \hat{\p}_0^{\tr}\E_1\p_0  & 0\\[.5em]
   0 & 0
 \end{pmatrix}\begin{pmatrix}
 \xi_p\\[.5em]
 \xi_q
 \end{pmatrix}'= \begin{pmatrix}
   \hat{\p}_0^{\tr}\A_0\p_0  & 0\\[.5em]
     \hat{\q}_0^{\tr}\A_0\p_0  & - \hat{\q}_0^{\tr}\E_1\q_0 
 \end{pmatrix}\begin{pmatrix}
 \xi_p\\[.5em]
 \xi_q
 \end{pmatrix}\\+\begin{pmatrix}
 \hat{\p}_0^{\tr} \f(\E_1\p_0\xi_p)\\[.5em]
  \hat{\q}_0^{\tr}\f(\E_1\p_0\xi_p)
 \end{pmatrix}
 +\begin{pmatrix}
  \hat{\p}_0^{\tr} \B\\[.5em]
  \hat{\q}_0^{\tr} \B
 \end{pmatrix}\u. \label{Eqn:13}
\end{multline}
The  system   \eqref{Eqn:13} can be \colb{r}educed to a nonlinear   decoupled system given by
\begin{subequations}
\label{eqn:decp111}
\begin{align}
{\bf E}_p\xi'_p&={\bf A}_p\xi_p+\f_p(\xi_p)+{\bf B}_p\bm{u},\quad \xi_p(0)=\p_0^{*\tr}\bm{x}(0),\label{eqn:decp11a}\\
{\bf E}_q\xi_q&={\bf A}_q\xi_p+\f_q(\xi_p)+{\bf B}_q\bm{u},\label{eqn:decp11b}\\
 \y&= \C_p\xi_p+\C_q\xi_q,
\end{align}
\end{subequations}
where 
\begin{align*}
 {\bf E}_p&=\hat{\p}_0^{\tr}{\bf E}_0\p_0\in \RR^{n_p\times n_p}, \, {\bf A}_p=\hat{\p}_0^{\tr}{\bf A}_0\p_0\in \RR^{n_p\times n_p},\, {\bf B}_p=\hat{\p}_0^{\tr}{\bf B}\in \RR^{n_p \times m}, \\
 {\bf E}_q&=-\hat{\q}_0^{\tr}{\bf A}_0\q_0\in \RR^{n_q \times n_q}, \, {\bf A}_q=\hat{\q}_0^{\tr}{\bf A}_0\p_0\in \RR^{n_p \times n_q}, \, {\bf B}_q=\hat{\q}_0^{\tr}\B\in \RR^{n_q \times m}.
\end{align*}
 The nonlinear terms  are defined as: $\f_p(\xi_p)=\hat{\p}_0^{\tr}\tilde{\f}(\xi_p) \in \RR^{n_p}, \, \f_q(\xi_p)=\hat{\q}_0^{\tr}\tilde{\f}(\xi_p)\in \RR^{n_q}$ where  $\tilde{\f}(\xi_p)=\f(\E_1\p_0\xi_p)\in \RR^n.$
 We note  that matrices $\E_p$ and $\E_q$ are always  nonsingular.
   We  can observe that  \eqref{eqn:decp111}   does not involve any matrix inversions.  It  is an implicit version of  the decoupled system \eqref{eqn:303} and   their output  solutions must coincide.
   However, in practice  it  is computationally cheaper to construct the coefficients of  \eqref{eqn:decp111}   than  \colp{those in} \eqref{eqn:303}.
Both decoupled systems preserve the dimension and the stability of the  nonlinear DAE \eqref{Eqn:1}. 
If \eqref{Eqn:1} is of tractability index $1$, then it can be  automatically decoupled  into either  \eqref{eqn:decp111} or \eqref{eqn:303}. Thus, instead of simulating 
\eqref{Eqn:1}, we can simulate its equivalent nonlinear decoupled system \eqref{eqn:decp111}  easily using standard numerical integration and  solvers.
 Decoupled systems \eqref{eqn:303} and \eqref{eqn:decp111} can be constructed in  efficient way by   employing  the sparse LU
decomposition-based routine, called LUQ, see \cite{Zhang}, to construct the projectors and their respective bases.
In the next section, we discuss how to apply MOR  \colp{to} \eqref{eqn:decp111}.
\section{Index-aware  MOR for nonlinear DAEs}
\label{sec:4}
Here, we consider  \colp{the}  equivalent nonlinear decoupled system \eqref{eqn:decp111} \colp{corresponding to the} nonlinear DAE \eqref{Eqn:1}\colp{,} but the same strategy can \colp{be} applied \colp{to} \eqref{eqn:303}.  
Given such a nonlinear decoupled system\colp{,} our goal \colp{is} to  reduce  the order of  differential  and algebraic parts separately.
\subsection{MOR for  the nonlinear differential subsystem}
We consider the nonlinear differential subsystem of  the nonlinear decoupled system \eqref{eqn:decp111} given by 
\begin{subequations}
\label{eqn:ode}
\begin{align}
{\bf E}_p\xi'_p&={\bf A}_p\xi_p+\f_p(\xi_p)+{\bf B}_p\bm{u},\quad \xi_p(0)=\p_0^{*\tr}\bm{x}(0),\label{eqn:odea}\\
 \y_p&= \C_p\xi_p,
\end{align}
\end{subequations}
where $\y_p\in \RR^{\ell\times n_p}$ is the output solution of the differential part. Our goal is reduction by projection of  system \eqref{eqn:ode}.  This means we want to find a linear subspace in which  the solution
trajectory lies approximately. This subspace is defined by \colp{its basis} matrix $\V_p\in \RR^{n_p\times r_p}$ where $r_p\ll n_p.$  We are interested in finding a solution  $\xi_{p_r}\in\mathbb{R}^{r_p}$  such that 
$\xi_p\approx \V_p \xi_{p_r}.$ We can then project  system \eqref{eqn:ode} onto that subspace by Galerkin projection resulting in the reduced differential subsystem
\begin{subequations}
\label{eqn:def}
\begin{align}
{\bf E}_{p_r}\xi'_{p_r}&={\bf A}_{p_r}\xi_{p_r}+\f_{p_r}(\xi_{p_r})+{\bf B}_{p_r}\bm{u}, \label{eqn:defa}\\
 \y_{p_r}&= \C_{p_r}\xi_{p_r},
\end{align}
\end{subequations}
where ${\bf E}_{p_r}=\V_p^{\tr} \E_p \V_p\in \mathbb{R}^{r_p\times r_p},\,{\bf A}_{p_r}=\V_p^{\tr} \E_p \V_p\in \mathbb{R}^{r_p\times r_p},$ $ {\bf B}_{p_r}=\V_p^{\tr} \B_p\in\mathbb{R}^{r_p\times m},$ \\$ \f_{p_r}(\xi_{p_r})=\V_p^{\tr} \f_p(\V_p \xi_{p_r})\in \mathbb{R}^{r_p}$ and $
\C_{p_r}=  \C_p\V_p\in \mathbb{R}^{\ell\times r_p}.$  Projection matrix $\V_p$ can be computed using standard MOR techniques for nonlinear systems such as POD \cite{morBenGW15}. However, if we employ POD by using 
 \eqref{eqn:defa} to compute the  
 snapshots, the nonlinearity $\f_p(\V_p \xi_{p_r})$ \colg{requires computation of $\f_p(\V_p \xi_{p_r})$ which has a complexity in the system dimension.}   \colp{Therefore}, we use  discrete empirical interpolation method (DEIM)  to create a truly low-dimensional function approximating  $\V_p^{\tr} \f_p(\V_p \xi_{p_r}).$  The DEIM algorithm creates matrices $\U_p,\W_p$ such that 
$$\V_p^{\tr} \f_p(\V_p \xi_{p_r})\approx \V_p^{\tr} \U_p(\W_p^{\tr}\U_p)^{-1}\W_p^\tr \f_p(\V_p \xi_{p_r}).$$
Here $\U_p\in \RR^{n_p\times m_p}$  is  orthonormal
 and the matrix  $\W_p\in \RR^{n_p\times m_p}$ is a picking matrix, where each row
has exactly one nonzero entry which is $1.$  This means that
$\W_p^\tr \f_p$ picks $m_p$ functions from the vector of functions $\f_p$. Here
we have to make sure to pick $m_p$ appropriately, in order to make
$\W_p^\tr \f_p(\V_p \xi_{p_r})$ truly low-dimensional, see \cite{morGruJ15}.
\subsection{Reduction of  algebraic subsystem}
  After reducing the differential subsystem using\colp{,} \colg{for example}\colp{,}  POD\colp{,}  the nonlinear term in the algebraic subsystem \eqref{eqn:decp11b} is also affected leading to 
\begin{subequations}
\label{eqn:alg}
\begin{align}
{\bf E}_{q}\xi_{q}&\approx {\bf A}_{q}\V_p \xi_{p_r}+\f_q(\V_p \xi_{p_r})+{\bf B}_{q}\bm{u},\label{eqn:alga}\\
 \y_q& \approx \C_{q}\xi_{q},
\end{align}
\end{subequations}
where $\y_q\in \RR^{\ell\times n_q}$ is the output solution of the algebraic part after reducing the differential subsystem.
Here, we intend to reduce the size of the algebraic variables $\xi_q$ by constructing another  matrix $\V_q\in \RR^{n_q\times r_q}$
where $r_q\ll n_q. $ Th\colp{at is}, we  replace \eqref{eqn:alg} by a reduced algebraic subsystem given by 
\begin{subequations}
\label{eqn:algr}
\begin{align}
{\bf E}_{q_r}\xi_{q_r}&= {\bf A}_{q_r}\xi_p+\f_{q_r}(\xi_{p_r})+{\bf B}_{q_r}\bm{u},\label{eqn:algr}\\
 \y_{q_r}&= \C_{q_r}\xi_{q_r},
\end{align}
\end{subequations}
where  ${\bf E}_{q_r}= \V_q^\tr \E_q\V_q\in \RR^{r_q\times r_q}, \, {\bf A}_{q_r}= \V_q^{\tr} {\bf A}_{q}\V_p \in \RR^{r_q\times r_p}, \, {\bf B}_{q_r}\in \RR^{r_q\times m},$\\ $ \C_{q_r}= \C_{q} \V_q \in \RR^{\ell\times r_q}$
and $\f_{q_r}(\xi_{p_r})= \V_q^{\tr}\f_q(\V_p \xi_{p_r}) \in \mathbb{R}^{r_q}.$
Reduction  matrix $\V_q$ can also be computed using the POD by  taking the algebraic solutions   of \eqref{eqn:decp11b}  \colp{ obtained from the snapshots of \eqref{eqn:decp11b} } as snapshots. Also here, the nonlinearity $\f_q(\V_p \xi_{p_r})$ has to be evaluated 
completely, even though we reduce the algebraic system size. Hence,  we also need to use  the DEIM to create a truly low-dimensional function approximating  $\V_q^{\tr}\f_q(\V_p \xi_{p_r}).$
The DEIM algorithm creates matrices $\U_q,\W_q$ such that 
$$\V_q^{\tr} \f_q(\V_q \xi_{p_r})\approx \V_q^{\tr} \U_q(\W_q^{\tr}\U_q)^{-1}\W_q^\tr \f_q(\V_p \xi_{p_r}),$$
where  $\U_q\in \RR^{n_q\times m_q}$ and $\W_q\in \RR^{n_q\times m_q}$ is a picking matrix. 
Combining  \eqref{eqn:def} and \eqref{eqn:algr}, we obtain an index-aware reduced order model (I-ROM) of \eqref{Eqn:1} given by
\begin{equation}
\label{Eqn:IROM}
 \begin{aligned}
  {\bf E}_{p_r}\xi'_{p_r}&={\bf A}_{p_r}\xi_{p_r}+\f_{p_r}(\xi_{p_r})+{\bf B}_{p_r}\bm{u}, \quad \xi_{p_r}(0)=\xi_{p_{r_0}},\\
  {\bf E}_{q_r}\xi_{q_r}&= {\bf A}_{q_r}\xi_p+\f_{q_r}(\xi_{p_r})+{\bf B}_{q_r}\bm{u},\\
 \y_{r}&= \C_{p_r}\xi_{p_r} + \C_{q_r}\xi_{q_r},
 \end{aligned}
 \end{equation}
where the reduced dimension is given by $r=r_p+r_q\ll n.$ Thus, we replace   \eqref{Eqn:1} with \eqref{Eqn:IROM}  instead of  \eqref{Eqn:Red}.

  
\section{Nonlinear DAEs arising from gas networks}
\label{sec:5}
In this section, we apply the  implicit decoupling strategy  proposed in Subsection \ref{sub:sec1}  \colp{to}  nonlinear DAEs arising from gas flow in pipeline  networks.

\subsection{Index reduction of DAEs arising from    {g}as  networks}
We consider   a spatial discretization approach of \colr{ one dimensional isothermal Euler equations arising from  gas flow pipe network\colp{s} }  proposed  in  \cite{morGruJHetal14,morGruHR16}\colp{,} leading to a nonlinear DAE given by 
\begin{subequations}
\label{Eqn:55}
 \begin{align}
  \vert \mathcal{A}_S^\tr\vert \partial_t \p_s +\vert \mathcal{A}_0^\tr\vert \partial_t \p_d  &=-\M_L^{-1}\q_{-},\\
  \partial_t \q_+&=\M_{A}(\mathcal{A}_S^\tr\p_s+\mathcal{A}_0^{\tr}\p_d)+\g(\q_+,\p_s,\p_d),\\
  0&= \mathcal{A}_0\q_+ +\vert \mathcal{A}_0\vert \q_{-}-\mathbf{B}_d\d(t),\\
  0&=\p_s-\s(t).
 \end{align}
\end{subequations}
The unknowns are  described {by} the pressure at the supply nodes $\mathbf{p}_s\in\mathbb{R}^{n_s}$, the pressure at all other nodes $\mathbf{p}_d\in\mathbb{R}^{n_d+n_0},$ the difference of flux over a pipe segment $\mathbf{q}_-\in\mathbb{R}^{n_E}$ and the average of the mass flux over a pipe segment $\mathbf{q}_+\in\mathbb{R}^{n_E}$, modelled over a graph with $n_E$ edge segments, that \colp{correspond to} the size of the discretization, $n_s$ supply nodes, $n_d$ demand nodes and $n_0$ interior nodes. 
\colb{The diagonal matrices } $\mathbf{M}_L\in \mathbb{R}^{n_E\times n_E}$ and $\mathbf{M}_A\in \mathbb{R}^{n_E\times n_E}$   encode parameters such as length, radius of the pipe segments as well as constants coming from the gas equation. The matrix $\mathbf{B}_d\in\mathbb{R}^{(n_d+n_0)\times n_d}$ is a matrix of ones and zeros making sure that the demand of the demand node is put at the right place in the mass flux equation.  
The matrix $\mathcal{A}_0 \in \mathbb{R}^{n_d\times n_{E}}$ is extracted from  the  incidence matrix  of the graph representing the refined gas transportation network and removing the rows corresponding to the supply nodes, while $\mathcal{A}_S \in \mathbb{R}^{n_s \times n_E}$ is 
the matrix extracted from  the incidence matrix by only taking rows corresponding  to the supply nodes.
 $\vert \mathcal{A}_0 \vert $ and $\vert \mathcal{A}_S \vert $ are the incidence matrices of the undirected graph defined as  the component-wise absolute values of the incidence matrices 
of the directed graph, see \cite{morGruHR16}.
  The input functions $\d(t)=(\ldots,d_i(t),\ldots)^\tr \in \mathbb{R}^{m_d} $ and $\s(t)=(\ldots,s_i(t),\ldots)^\tr \in \mathbb{R}^{m_s}$ are vectors for flux (mass flow) at demand nodes  and  pressure at supply  nodes, respectively. 
  The nonlinear term $\g(  \q_+,
 \p_d,
 \p_s)=(\ldots,g_k(  \q_+,
 \p_d,
 \p_s),\ldots)^\tr\in \mathbb{R}^{n_E},$  is the  vector involving friction and gravitation effects {with} 
 \begin{align}
  g_k(  \q_+,
 \p_d,
 \p_s)&=-\frac{gA_k}{2\gamma_0}\psi_k(\p_d,\p_s)\frac{\Delta h_k}{L_k}-\frac{\lambda_k\gamma_0}{4 D_k A_k}\frac{\q_{+}^k\vert\q_{+}^k \vert }{\psi_k (\p_d,\p_s)}, \label{Eqn:gx}
 \end{align}
where $\psi_k (\p_d,\p_s)$ is the k-th entry of the vector-valued function: 
$$\psi(\p_d,\p_s)= \vert \mathcal{A}_S^\tr\vert  \p_s +\vert \mathcal{A}_0^\tr\vert \p_d \in \mathbb{R}^{n_E}.$$
The scalars $\lambda_k, D_k, L_k$ and $A_k$ denote friction,  diameter, length and area of the pipe's $k$-th  segment.  The scalar $\Delta h_k$ denotes the height difference of the pipe segment.  
These  scalar parameters in  the system and those defined earlier  are 
known at least within some range of uncertainty. 
 System  \eqref{Eqn:55} can be rewritten in the form \eqref{Eqn:1} leading to  a system of  nonlinear DAEs  with dimension $n=2n_E+n_d+n_0+n_s.$ 
The desire{d} outputs  in $\mathbb{R}^{n_s+n_d}$ can be obtained using the output equation
\begin{align}
 \mathbf{y}&={\begin{pmatrix} 
              \mathbf{y}_q\\
              \mathbf{y}_p
             \end{pmatrix}}=
\begin{pmatrix}
     0 & \vert \mathcal{{A}}_S\vert & 0 &0\\
     0& 0& \mathbf{B}_d^T&0
    \end{pmatrix}\begin{pmatrix}
  \mathbf{q}_-\\
   \mathbf{q}_+\\
 \mathbf{p}_d\\
 \mathbf{p}_s\\
 \end{pmatrix},\label{Eqn:output}
\end{align}
where $\mathbf{y}_q=\vert \mathcal{{A}}_S\vert\mathbf{q}_+$  is the mass flow at the supply nodes  and $\mathbf{y}_p=\mathbf{B}_d^T\mathbf{p}_d$ is the pressure at demand nodes.  We can observe that the 
 initial condition has to be consistent with the hidden constraints  in \eqref{Eqn:55}.
Efficient simulation of \eqref{Eqn:55} has numerical integration challenges  since the solutions of hyperbolic balance laws can blow-up in finite time, due to  both the stiffness and index problem.
In \cite{morGruJHetal14}, an index reduction strategy was proposed to eliminate the index problem. This was done by reformulating   \eqref{Eqn:55}  into  an implicit  nonlinear ODE given by
\begin{multline}
 \begin{pmatrix}
   \vert \mathcal{A}_0\vert \mathbf{M}_L \vert \mathcal{A}_0^\mathrm{T}\vert && 0\\
   0 && \mathrm{I}
 \end{pmatrix}\begin{pmatrix}
 {\partial_t} \mathbf{p}_d\\
   {\partial_t} \mathbf{q}_+
 \end{pmatrix}
 =\begin{pmatrix}
    0 && \mathcal{A}_0\\
    \mathbf{M}_A \mathcal{A}_0^{T} && 0
   \end{pmatrix}\begin{pmatrix}
  \mathbf{p}_d\\
 \mathbf{q}_+
 \end{pmatrix}+\begin{pmatrix}
  \vert \mathcal{A}_0\vert \mathbf{M}_L \vert \mathcal{A}_S^\mathrm{T}\vert {\partial_t} \mathbf{s}(t)\\
  \mathbf{g}(\mathbf{q}_+,\mathbf{s}(t),\mathbf{p}_d)
 \end{pmatrix}\\+\begin{pmatrix}
 0 && -\mathbf{B}_d\\
    \mathbf{M}_A \mathcal{A}_S^{T} && 0
 \end{pmatrix}\begin{pmatrix}
 \mathbf{s}(t)\\
 \mathbf{d}(t)
 \end{pmatrix}.\label{ode}
\end{multline}
 Since from \eqref{Eqn:output} we are  just interested in the  solutions of   $  \mathbf{q}_+$ and $ \mathbf{p}_q\colp{,}$    the dimension of the nonlinear  DAE \eqref{Eqn:55} can  be reduced to $\tilde{n}= n_d+n_0+n_E$ with output equation
  $$\mathbf{y}={\begin{pmatrix} 
              \mathbf{y}_q\\
              \mathbf{y}_p
             \end{pmatrix}}=
\begin{pmatrix}
      0& \vert \mathcal{{A}}_S\vert  \\
      \mathbf{B}_d^T&0 
    \end{pmatrix}\begin{pmatrix}
     \mathbf{p}_d\\
   \mathbf{q}_+
 \end{pmatrix}.$$
The generated ODE   can be  reduced further using  
standard MOR methods for nonlinear systems, such as POD, POD-DEIM, etc,  applied \colp{to}   \eqref{ode}, see \cite{morGruHR16}. 
However, the index reduction  approach presented  depends on the spatial discretization approach   used.
\colr{
In the next section, we propose an alternative  model which preserves the DAE structure independent of the spatial discretization method. }



\subsection{Decoupled  model of   {g}as {t}ransport networks}
Here, we discuss the decoupling analysis of nonlinear DAE \eqref{Eqn:55} arising from the gas transportation networks. \colr{As a result, we present an alternative model to the ODE model \eqref{ode} proposed in \cite{morGruHR16}.}
We can observe that  \eqref{Eqn:55} can \colp{be} re\colp{-}written into the form \eqref{Eqn:1} where

\begin{align}
\nonumber
\E&=\begin{pmatrix}
  0& 0&\vert \mathcal{A}_0^\tr\vert& \vert\mathcal{A}_S^\tr\vert\\
   0&\mathrm{I}&0&0\\
 0 & 0& 0& 0\\
   0& 0& 0 & 0
 \end{pmatrix},\quad   \A=\begin{pmatrix}
-\M_L^{-1} &0&0&0\\
   0&0&\M_{A}\mathcal{A}_0^\tr& \M_{A}\mathcal{A}_S^\tr\\
   \vert \mathcal{A}_0\vert &  \mathcal{A}_0 & 0& 0 \\
  0& 0&0&  \mathrm{I}
 \end{pmatrix},\\
 \B&=-\begin{pmatrix}
 0& 0\\
 0&  0\\
 0&\mathbf{B}_d\\
 \mathrm{I}&0
 \end{pmatrix},\quad  \C=\begin{pmatrix}
     0 & \vert \mathcal{A}_S\vert & 0 &0\\
     0& 0& \mathbf{B}_d^\tr&0
    \end{pmatrix},\quad  \f(\E\x)= \begin{pmatrix}
 0\\\tilde{\g}(\E\x)\\
 0\\
0
 \end{pmatrix}. \label{Eqn:60}
\end{align}
The unknown vector $\x$  and input vector $\u$ are given by 
  $
 \x=\begin{pmatrix}
  \q_-^\tr &
   \q_+^\tr&
 \p_d^\tr&
 \p_s^\tr
 \end{pmatrix}^\tr$ and $
    \u= \begin{pmatrix}
 \s(t)^\tr &
 \d(t)^\tr
 \end{pmatrix}^\tr,
 $ respectively.\\
$\tilde{\g}(\E\x)=\colb{\tilde{\g}(\psi_k(\p_d,\p_s), \q_+)=}(\ldots,\tilde{g}_k(\psi_k(\p_d,\p_s), \q_+),\ldots)^\tr,$ where $$\tilde{g}_k(\psi_k(\p_d,\p_s), \q_+)=-\frac{gA_k}{2\mu_0}\psi_k(\p_d,\p_s)\frac{\Delta h_k}{L_k}-\frac{\lambda_k\mu_0}{4 D_k A_k}\frac{\q_{+}\vert\q_{+} \vert }{\psi_k (\p_d,\p_s)}.$$
Since the  gas transport model can be rewritten in the form \eqref{Eqn:1},   we can decoupled it into either the form  \eqref{eqn:303} or   \eqref{eqn:decp111}.
    In our discussion, we shall use the implicit decoupling strategy proposed in subsection \ref{sub:sec1}
     leading to an implicit decoupled system \eqref{eqn:decp111}.
 For convenience, we can \colp{p}artition  \eqref{Eqn:60}  into  a   block form  leading to
 \begin{subequations}
 \label{Eqn:9}
  \begin{align}
  \begin{pmatrix}
   0& 0& \E_{13}\\
   0&\mathrm{I}& 0\\
   0& 0&0
  \end{pmatrix}\begin{pmatrix}
  \x_1\\
  \x_2\\
  \x_3
  \end{pmatrix}'&=\begin{pmatrix}
   \A_{11}&  0 & 0\\
   0&0&\A_{23}\\
   \A_{31}&\A_{32}&\A_{33}
  \end{pmatrix}\begin{pmatrix}
  \x_1\\
  \x_2\\
  \x_3
  \end{pmatrix}
  +
  \begin{pmatrix}
  0\\
 \tilde{\g}(\E\x)\\
  0
  \end{pmatrix}
  +\begin{pmatrix}
  0\\
  0\\
  \B_3
  \end{pmatrix}\begin{pmatrix}
 \s(t)\\
 \d(t)
 \end{pmatrix}, \label{Eqn:parta}\\
  \y&=\begin{pmatrix}
       0& \C_2 & \C_3
      \end{pmatrix}\begin{pmatrix}
  \x_1\\
  \x_2\\
  \x_3
  \end{pmatrix},
 \end{align}
  \end{subequations}
  where  $\E_{13}=\begin{pmatrix}
                \vert \mathcal{A}_0^\tr\vert& \vert \mathcal{A}_S^\tr\vert
               \end{pmatrix}\in \mathbb{R}^{n_E\times n_v},\quad \A_{11}=
                -\M_L^{-1}\in \mathbb{R}^{n_E\times n_E},$ \\$
               \A_{23}=\begin{pmatrix}
               \M_{A}\mathcal{A}_0^\tr& \M_{A}\mathcal{A}_S^\tr
               \end{pmatrix} \in \mathbb{R}^{n_E\times n_v},\quad   \A_{31}=\begin{pmatrix}
   \vert \mathcal{A}_0\vert \\
   0
  \end{pmatrix}\in \mathbb{R}^{n_v \times n_E},$ \\$   \A_{32}=\begin{pmatrix}
     \mathcal{A}_0\\
   0
  \end{pmatrix}\in \mathbb{R}^{n_v\times n_E},\,
  \A_{33}=\begin{pmatrix}
   0 &  0\\
   0& \mathrm{I}
  \end{pmatrix}\in \mathbb{R}^{n_v\times n_v}, \,\B_3=-\begin{pmatrix}
   0&\colr{\B_d}\\
 \mathrm{I}&0
  \end{pmatrix} \in \mathbb{R}^{n_v\times m},$\\ $\C_2=\begin{pmatrix}
      \vert \mathcal{A}_S\vert \\
      0
    \end{pmatrix}\in \mathbb{R}^{\ell\times n_E},\, \C_3=\begin{pmatrix}
 0 &0\\
    \mathbf{B}_d^T&0
    \end{pmatrix}\in \mathbb{R}^{\ell\times n_v},\, \x_1=\q_{-} \in \mathbb{R}^{n_E},$\\ $  \x_2=\q_{+}\in \mathbb{R}^{n_E},$ $ \x_3=\begin{pmatrix}
                                   \p_d\\
                                   \p_s
                                  \end{pmatrix}\in \mathbb{R}^{ n_v},\,  
$ $n_v=n_d+n_s+n_0.$
The nonlinear term is defined as   $$\tilde{\g}(\E\x)=\tilde{\g}(\E_{13}\x_3,\x_2,0,0)=\tilde{ \tilde{\g}}(\x_3,\x_2)=(\ldots,  \tilde{\tilde{g}}_k(\x_3^k,\x_2^k),\ldots)^\tr\in \mathbb{R}^{n_E},$$ {with} 
 \begin{align}
  \tilde{\tilde{g}}_k(\x_3^k,\x_2^k)&=-\frac{gA_k}{2\gamma_0}\E_{13}\x_3^k\frac{\Delta h_k}{L_k}-\frac{\lambda_k\gamma_0}{4 D_k A_k}\frac{\x_2^k\vert\x_2^k \vert }{\E_{13}\x_3^k}. 
 \end{align} 
In order to decouple \eqref{Eqn:9}, we need to   first find  the tractability index of \eqref{Eqn:9}  using Definition \ref{defn1}.
Setting 
\begin{align}
 \label{eqn:E0}
 \E_0&= \begin{pmatrix}
   0& 0& \E_{13}\\
   0&\mathrm{I}& 0\\
   0& 0&0
  \end{pmatrix} \quad \mbox{and}\quad \A_0=\begin{pmatrix}
   \A_{11}&  0 & 0\\
   0&0 &\A_{23}\\
   \A_{31}&\A_{32}&\A_{33}
  \end{pmatrix},
\end{align}
 \colb{w}e can then  construct    projectors
  \begin{align}
  \label{Eqn:Q0}
   \Q_0&=\begin{pmatrix}
    \mathrm{I}& 0& 0\\
   0&0& 0\\
   0& 0& \Q
  \end{pmatrix}\in \mathbb{R}^{n\times n} \quad \mbox{and}\quad   \P_0=\mathrm{I}-\Q_0=\begin{pmatrix}
    0&0& 0\\
   0&  \mathrm{I}&  0\\
   0& 0& \P
  \end{pmatrix}\in \mathbb{R}^{n\times n},
  \end{align}
 such that $\E_0\Q_0=0,$ \colg{meaning $\E_{13}\Q=0$}, or
 $\Q\in \mathbb{R}^{ n_v \times n_v}$ is  the projector onto the nullspace of $\E_{13}$ and $\P\in \mathbb{R}^{ n_v \times n_v}$ is its complementary projector.
 Substituting  the above matrices and projectors into \eqref{proj} leads \colp{to}
  $$\E_1=\E_0-\A_0\Q_0=\begin{pmatrix}
   -\A_{11}&  0 & \E_{13}\\
   0&\mathrm{I}&-\A_{23}\Q\\
   -\A_{31}&0&-\A_{33}\Q\end{pmatrix}.
 $$
 \colp{If} $\E_1$ is nonsingular,    the   DAE \eqref{Eqn:9}  is of tractability index $1$.
Next, we construct the values of the matrix coefficients of 
  \eqref{eqn:decp111} as follows.    \colb{Let $n_p=\mathrm{rank(\E_0)}$ and $n_q=n-n_p.$}
\colb{Then, t}he columns of the  matrices 
\begin{align}
 \label{eqn:q0}
 \q_0&=\begin{pmatrix}
        \mathrm{I}&  0\\
   0& 0\\
   0&  \q
      \end{pmatrix}\in \RR^{n\times n_q}\quad  \mbox{and}\quad   \p_0=\begin{pmatrix}
        0& 0\\
    \mathrm{I}&  0\\
    0& \p
      \end{pmatrix}\in \RR^{n\times n_p}
\end{align}
  are linearly independent and span the column spaces  of $\Q_0$  and $\P_0$ in \eqref{Eqn:Q0}, respectively.   
     The left inverse of column matrices  $\p_0$  and $\q_0$ are given by 
     \begin{align}
     \label{eqn:q1}
      \q_0^{*\tr}&=\begin{pmatrix}
        \mathrm{I}&  0 & 0\\
   0& 
   0&  \q^{*\tr}
      \end{pmatrix}\in \RR^{n_q\times n}\quad  \mbox{and}\quad   \p_0^{*\tr}=\begin{pmatrix}
        0&  \mathrm{I}&  0 \\
   0& 
   0&  \p^{*\tr}
      \end{pmatrix}\in \RR^{ n_p \times n},
     \end{align}
      respectively, where 
      $\q^{*\tr}$ and $\p^{*\tr}$ are  the  left inverse\colp{s} of column matrices  $\q$ and $\p,$ respectively.    Let $k_q$ be the dimension of  the nullspace of $\E_{13},$ and 
      $k_p=n_v-k_q.$ The columns of   $\q\in \mathbb{R}^{n_v\times k_q}
$ and $\p\in \mathbb{R}^{n_v\times k_p}$  are linearly independent and span the column spaces  of  $\Q$  and $\P$ in \eqref{Eqn:Q0}, respectively.
Finally,  \colr{column}
matrices \colr{$\hat{\p}_0\in \mathbb{R}^{n\times n_p}$ and $\hat{\q}_0\in \mathbb{R}^{n\times n_q}$}  can be constructed such that   their  columns are linearly independent \colr{and } span the null spaces of  the matrices   $\q_0^{\tr}\A_0^{\tr}\in \RR^{n_q\times n}$  and    $\E_0^{\tr}\in \RR^{n\times n},$ respectively.    
The differential  and algebraic variables are given by 
   \begin{align}
       \label{eqn:q2}
       \xi_p= \p_0^{*\tr} \P_0\x=\begin{pmatrix}
                                                     \x_2\\
                                                     \p^{*\tr}\x_3
                                                    \end{pmatrix}\quad  \mbox{and}\quad   \xi_q =\q_0^{*\tr} \Q_0\x=\begin{pmatrix}
                                                     \x_1\\
                                                     \q^{*\tr}\x_3
                                                    \end{pmatrix},
   \end{align}
 respectively.   The nonlinear term is defined as 
 
 \begin{align}
 \label{eqn:f}
  \tilde{\f}(\xi_p)&= \begin{pmatrix}
  0\\
  \g_p(\xi_p)\\
  0
  \end{pmatrix}
 \end{align}
 where     
 $$\g_p(\xi_p)=\tilde{\g}(\E_1\p_0\xi_p)=\tilde{\g}(\E_{13}\p\xi_{p_2}, \xi_{p_1},0)=\tilde{\tilde{\g}}(\xi_{p_1},\xi_{p_2})=(\ldots,  \tilde{\tilde{g}}_k(\xi_{p_1}^k,\xi_{p_2}^k),\ldots)^\tr\in \mathbb{R}^{n_E},$$ {with} 
$$  \tilde{\tilde{g}}_k(\xi_{p_1}^k,\xi_{p_2}^k)=-\frac{gA_k}{2\gamma_0}\tilde{\E}_{13} \xi_{p_2}^k\frac{\Delta h_k}{L_k}-\frac{\lambda_k\gamma_0}{4 D_k A_k}\frac{\xi_{p_1}^k\vert\xi_{p_1}^k \vert }{\tilde{\E}_{13} \xi_{p_2}^k}\quad \mbox{and}\quad \tilde{\E}_{13}= \E_{13}\p.$$ 
This is due to the fact that $\E_{13}\x_3=\E_{13}\p\p^{*\tr}\x_3=\E_{13}\p\xi_{p_2}.$   It can  be proved that $\f_q(\xi_p)=\hat{\q}_0^{\tr}\tilde{\f}(\xi_p)=0\in \RR^{n_q}$ always due to the structure of the nonlinearity.
 Finally, substituting \eqref{eqn:E0}, \eqref{eqn:q0}-\eqref{eqn:q2} into \eqref{eqn:decp111} leads to an equivalent  nonlinear decoupled system of \eqref{Eqn:60}  given by 
\begin{subequations}
\label{eqn:decp111s}
\begin{align}
{\bf E}_p\xi'_p&={\bf A}_p\xi_p+\f_p(\xi_p)+{\bf B}_p\bm{u},\quad \xi_p(0)=\begin{pmatrix}
                                                     \x_2(0)\\
                                                     \p^{*\tr}\x_3(0)
                                                    \end{pmatrix},\label{eqn:decp11sa}\\
{\bf E}_q\xi_q&={\bf A}_q\xi_p+{\bf B}_q\bm{u},\label{eqn:decp11bs}\\
 \y&= \C_p\xi_p+\C_q\xi_q, \label{eqn:decpss}
\end{align}
\end{subequations}
where $\f_p(\xi_p)=\hat{\p}_0^{\tr}\tilde{\f}(\xi_p)\in \RR^{n_p}$  with  $\tilde{\f}(\xi_p)$ as defined in \eqref{eqn:f}. \colr{The system matrix coefficients are computed   as  defined  in \eqref{eqn:decp111}.}
We can observe  that the nonlinear  gas  transport network model has been decoupled into $n_p=n_E +k_p$ nonlinear differential equations, and $n_q=n_E+k_q$ algebraic equations.
This decoupled system preserves all the physical properties of  the DAE \eqref{Eqn:55} such as hyperbolicity.
Subsystem \eqref{eqn:decp11sa} can be simulated using  standard numerical integration, then algebraic solutions of \eqref{eqn:decp11bs} can be obtained by using numerical solvers after  post-processing. 
Hence, the desired  output data can be obtained through \eqref{eqn:decpss}.
We note \colp{that} the decoupling  enable\colp{s} us to treat DAEs \colp{like} ODEs\colp{.} \colp{H}owever\colp{,} the stiffness problem is inherited in the ODE subsystem \eqref{eqn:decp11sa}. 
In order to \colp{cope with} the stiffness problem, we can use  IMEX integration scheme \cite{morGruJ15} \colp{instead of} standard integration which
makes an efficient simulation of \eqref{eqn:decp11sa}  possible. We note that  the values of the matrix coefficients  of  \eqref{eqn:decp111s} can vary depending on the choices of projectors in \eqref{Eqn:Q0}\colp{,} but the solutions will always be the same. 
In practice, system \colp{ \eqref{eqn:decp111s}} can be constructed automatically following the implicit decoupling procedure in Subsection \ref{sub:sec1}. Numerical experiments show that \eqref{eqn:decp11sa} and \eqref{ode} have the same dimension for the case of index $1$ gas transportation
networks.

\section{Numerical experiments}
\label{sec:num}
In this section, we illustrate  the performance of the proposed decoupling and IMOR method  for nonlinear DAEs with a special nonlinear term $\f(\x)=\f(\E\x),$ where $\E$ is a singular matrix.
Such nonlinear DAEs  can  arise from  gas transportation  networks  as discussed in Section \ref{sec:5}. Here, we consider small to large  examples of gas transportation  networks  leading to nonlinear DAEs of tractability index $1.$
%
We compute the relative error in the  format 
 $ \mathrm{Re. error}=\Vert \mathbf{y}-\mathbf{y}_r\Vert_2/\Vert \mathbf{y}\Vert_2.  $
 The output error is defined as $\mathrm{max}( \mathrm{Re. error} (\mbox{pressure}), \mathrm{Re. error} (\mbox{mass flow}) ).$
 {Simulations  were done using   MATLAB\textregistered Version 2012b on a Unix desktop.}

 \subsection{Numerical  integration}
 We compare  the output solutions (mass flow at \colb{the} supply node and pressure at demand nodes) of different gas transportation models:  nonlinear  DAE model \eqref{Eqn:55},  nonlinear ODE  model \eqref{ode} and   nonlinear  decouple\colb{d} model \eqref{eqn:decp111s}.
\begin{example}
 \label{Examp:1}
 In this example, we consider small to medium gas pipeline  networks   from \cite{Kaarthik18,morGruHKetal13} with  steady pressure at \colp{the} supply  pressure node and steady mass flow at demand nodes.
 We are interested in the comparison of the pressure and mass flows of different models  of  each gas transportation network shown in Table \ref{Tab:banagaaya_5}.
     \begin{table}[!h]
 \footnotesize
  \caption{Comparison of gas transportation models }
  \label{Tab:banagaaya_5}
 \begin{center}
 \scalebox{1}{
\begin{tabular}{lllllll}
\hline 
\textbf{Nonl.   DAE } & \textbf{ Nonl.  ODE  } & \multicolumn{3}{l}{ \textbf{ Nonl.  Decoupled }}   & \textbf{Supply nodes} &\textbf{ Demand nodes}\\
$n$ & $\tilde{n}$ & $n_p$ & $n_q$ & $n_p+n_q$ & $m_s$ & $m_d$\\ \hline 
$4$ & $2$  & $2$ & $2$ &$4$& $1$ & $1$\\
$25$ & $16$  & $16$ & $9$ &$25$& $1$ & $2$\\
$55$ & $36$  & $36$ & $19$ &$55$& $1$ & $8$\\
$121$ & $80$  & $80$ & $41$ &$121$& $1$ & $24$\\ \hline 
\end{tabular}
}
\end{center}
 \end{table}

\ \\
  In Table \ref{Tab:banagaaya_5}, we can observe that  the index reduced   ODE model has the same dimension   as the differential 
 part of the nonlinear decoupled model. We use the implicit-Euler numerical integration scheme to solve the nonlinear DAE and ODE models with a fixed time step. For the nonlinear  decoupled model we use the  implicit-Euler numerical integration scheme 
 on the differential part and LU based numerical solver \colp{for} the algebraic part.
  Figures \ref{Fig_1:1}-\ref{Fig_1:4}, show the pressure at \colp{the} supply node,  mass flow at the first demand node, mass flow at the supply node and pressure at the first demand \colp{node} for each network presented in Table \ref{Tab:banagaaya_5}. 
 In Figure \ref{Fig_1:1}, we used  \colb{steady} pressure $\mathbf{s}(t)=650 \mathrm{bars}$ at the supply node and \colb{steady}  mass flow rate of  $\mathbf{d}(t)=100 \mathrm{Kg/s}$ at the demand node. 
 \begin{figure}[!h]
 \centering
  \includegraphics[width=\textwidth,height=.6\textwidth]{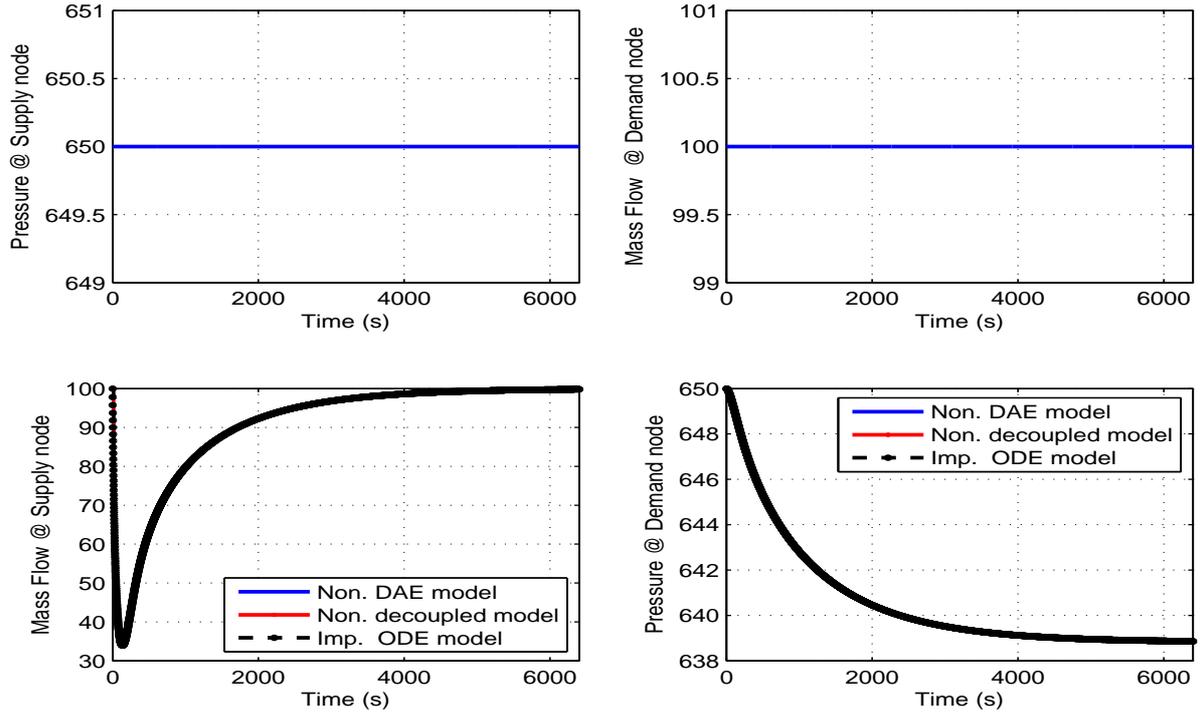} 
  \caption{Gas transportation network $(n=4, m_s=1, m_d=1)$}
  \label{Fig_1:1}
 \end{figure}
 
 \ \\
 In Figure \ref{Fig_1:2}, we used \colb{steady} pressure $\mathbf{s}(t)=700 \mathrm{bars}$ at the supply node and \colb{steady} mass flow rate of  $\mathbf{d}(t)=(60,30)^\tr$ at the demand nodes.

  \begin{figure}[!h]
   \centering
  \includegraphics[width=\textwidth,height=.6\textwidth]{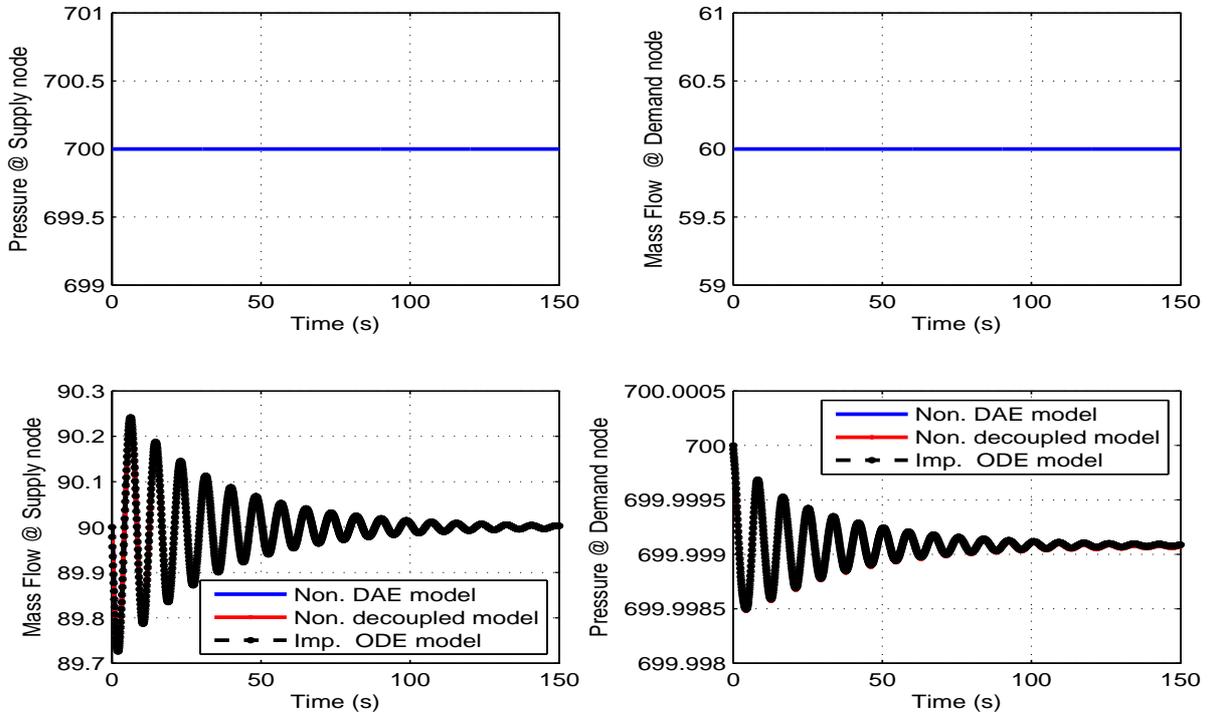}
   \caption{Gas transportation network $(n=25, m_s=1, m_d=2)$}
     \label{Fig_1:2}
 \end{figure}
  \begin{figure}[!h]
   \centering
  \includegraphics[width=\textwidth,height=.65\textwidth]{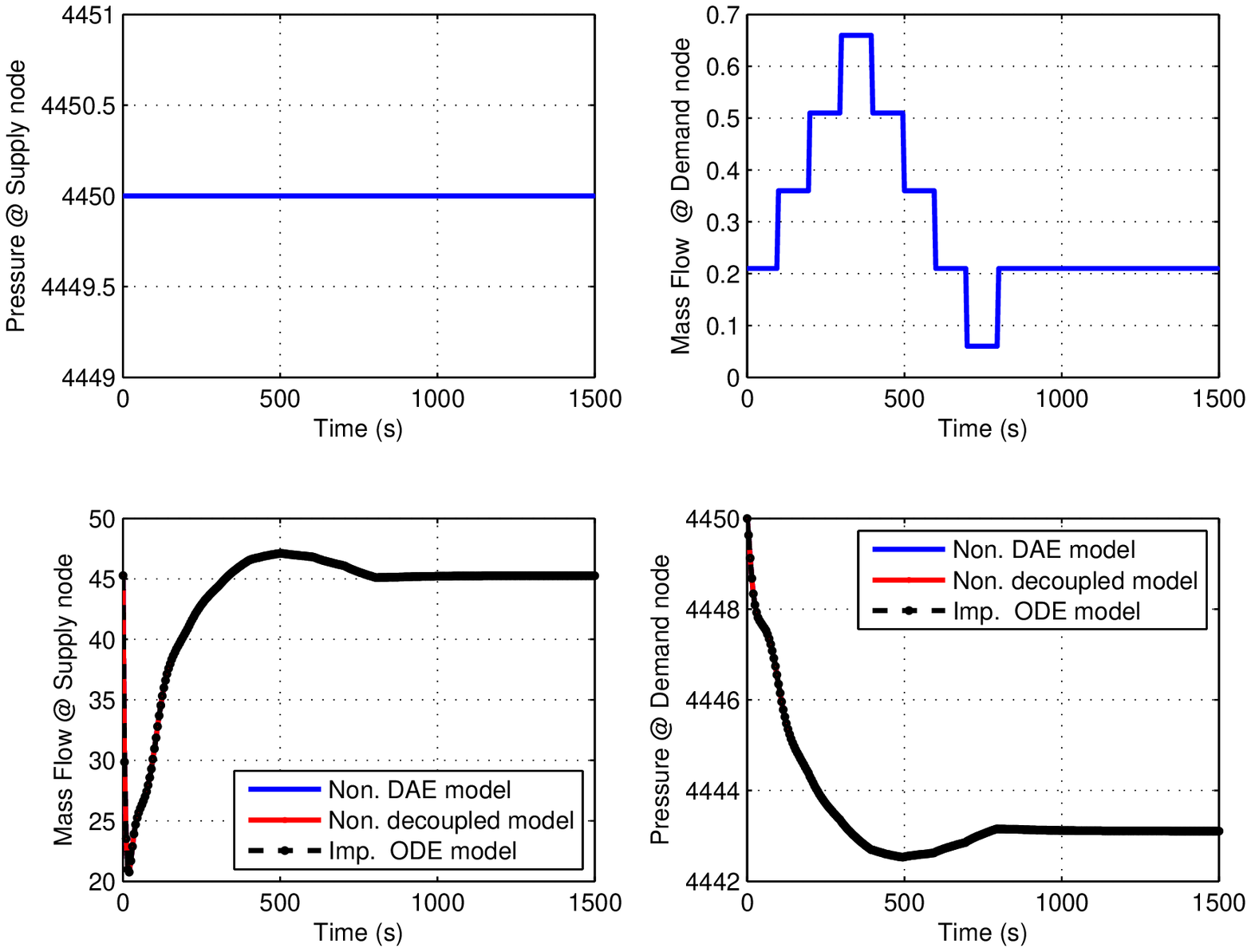}
   \caption{Gas transportation network $(n=55, m_s=1, m_d=8)$}
     \label{Fig_1:3}
 \end{figure}

 \begin{figure}[!h]
  \centering
  \includegraphics[width=\textwidth, height=.65\textwidth]{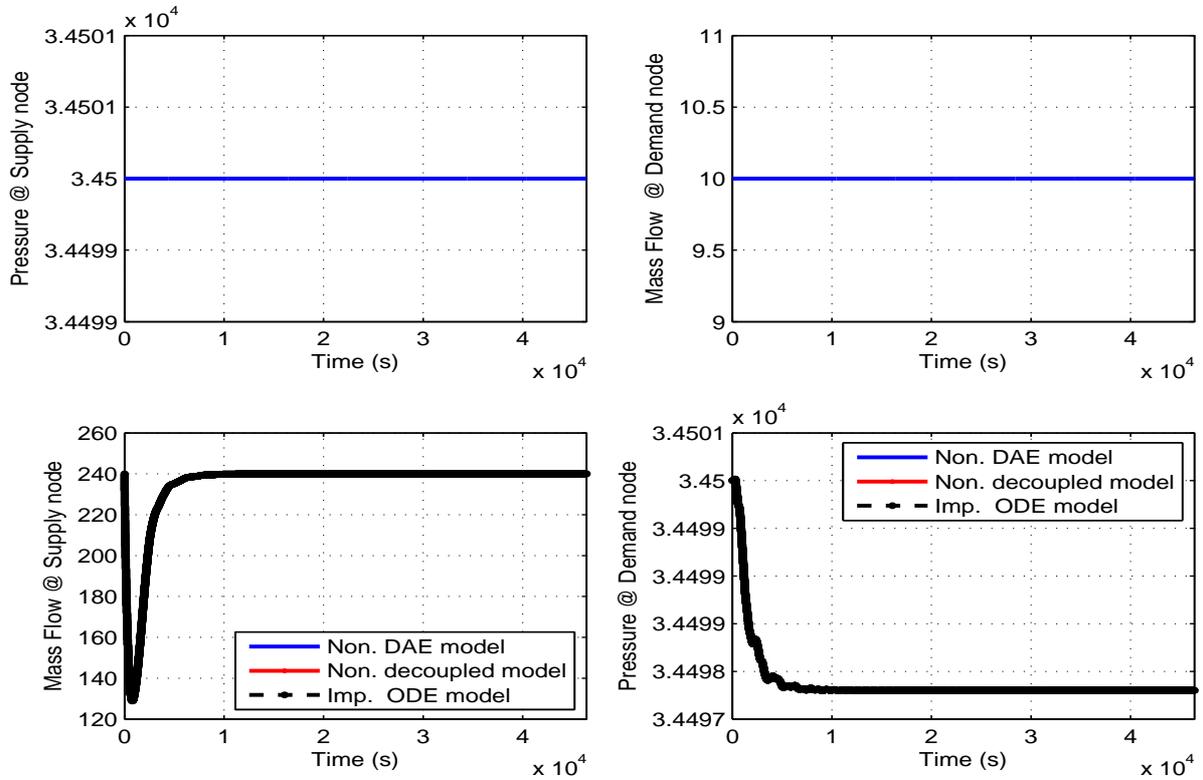}
   \caption{Gas transportation network $(n=121, m_s=1, m_d=24)$}
     \label{Fig_1:4}
 \end{figure}

  \ \\
  In Figure \ref{Fig_1:3}, we used \colb{steady} pressure $\mathbf{s}(t)=4.55\times 10^4 \mathrm{bars}$ at the supply node and \colb{steady} mass flow rate of     $$\mathbf{d}(t)=(2.1,348.6,2.2,28.3,18.1,10.4,28.5,14.5)^\tr$$ at the demand nodes.
   In Figure \ref{Fig_1:4}, we used \colb{steady} pressure $\mathbf{s}(t)=3.45\times 10^4 \mathrm{bars}$ at the supply node and \colb{steady} mass flow rate of  $\mathbf{d}(t)=10\times \mathrm{ones} (24,1)$ at the demand nodes. 
In all test cases, we can observe that all  models decay towards   steady  mass flow at the supply node and  steady pressure at the demand no\colb{d}es. 

\end{example}

\begin{example}
\label{Examp:1b}
In this  example, we are interested in comparing the pressure and   mass flow rate while applying   steady pressure  at supply node and  transient  mass flow rate  at \colp{the} demand \colp{node}. 
We consider a medium size   gas transport pipe network  with $200$ pipes, one supply node and one demand node   generated using the following data. 
The length, diameter and average roughness of each pipe are chosen \colp{as} constants  given by $18.15 \mathrm{m}, \, 1.422\mathrm{m}$ and $1.5\times 10^{-6}\mathrm{m},$ respectively. 
The gas composition through the network  is methane with  specific gas constant $518.26 \mathrm{J/KgK}$ at steady supply of  $84\mathrm{bar}$ and 
mass flow at demand as shown in the first row of  Figure  \ref{banagaaya_1b:1} in the time interval $t\in \left[0,1000s \right].$
\begin{figure}[!h]
  \centering
  \includegraphics[width=\textwidth,height=.65\textwidth]{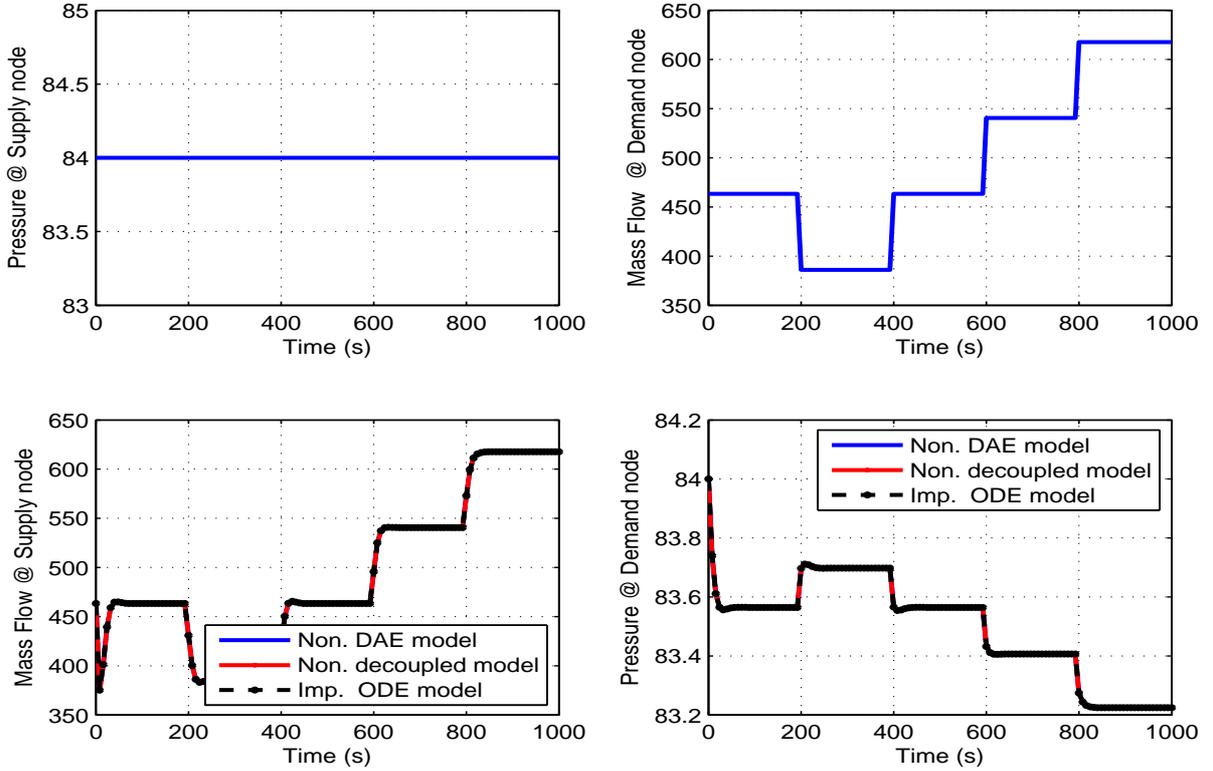}
  \caption{Comp\colb{a}rison of the output solutions }
  \label{banagaaya_1b:1}
 \end{figure}

 \ \\
This leads to a nonlinear DAE system of dimension $n=601$ which we  decoupled into  $n_p=400$ differential equations 
and $n_q=201$ algebraic equations.  For comparison, we  generated the ODE model  \eqref{ode}  leading   to an ODE model of dimension $400.$
In all models for integration, we use the implicit-Euler scheme \colp{with} the same step size of $8.$ 
In the second row of Figure \ref{banagaaya_1b:1}, we can observe that the \colr{  the pressure and mass flow coincide  with the  nonlinear DAE model  for both ODE model and the decoupled model. 
Using the solutions of the nonlinear DAE model as  reference,  the  solutions from the ODE model have  relative errors  of $2.4\times 10^{-6}$ and $5.2\times 10^{-8}$ in the pressure and mass flow, respectively, while 
the solutions from the decoupled model   have relative errors of $3.1\times 10^{-6}$ and $4.5\times 10^{-7},$ respectively. 
}  
 


 \end{example}
 
\begin{example}

\label{Examp:2}

  In this example, we consider a small size gas transport network  obtained  from \cite{morGruHKetal13}. It consists of $17$ nodes, $16$ pipes,
$1$ supply node and $8$ demand nodes.  Spatial discretization leads to  a nonlinear DAE  with 
 $$n=55,  m=\ell=9, m_s=1, m_d=8.$$ 
  We used steady pressure of  $\mathbf{s}(t)=4450 \mathrm{bars}$ at the supply node and mass flow rate of     $\mathbf{d}(t)=(0.21,34.86,0.22,2.83,1.81,1.04,2.85,1.45)^\tr$ at the demand nodes.
The  nonlinear  implicit ODE model  \eqref{ode} lead\colp{s} to a system of dimension $36$ 
while the  decoupled system  \eqref{eqn:decp111s} has   $n_p=36$ 
 differential equations 
and $n_q=19$ algebraic equations. 
We used the implicit\colp{-}Euler integration scheme to simulate {the} linear DAE  and implicit ODE models.   We also used the same method to simulate the ODE part and the LU method  for  {solving} the 
 algebraic part of the decoupled system. Using the same time steps  and time interval, we simulated all the models and some of the results are presented in Figure \ref{Banagaaya_fig3:1}. 
 In   Figure \ref{Banagaaya_fig3:1},   we only present  pressure and mass flow at \colp{the} supply node,  mass flow and pressure    at \colp{the} first  demand node.
 We can observe {that} the solutions of  the  nonlinear DAE \colr{model concides with both the ODE and decoupled models.} 
 \begin{figure}[h]
\includegraphics[width=\textwidth, height=.65\textwidth]{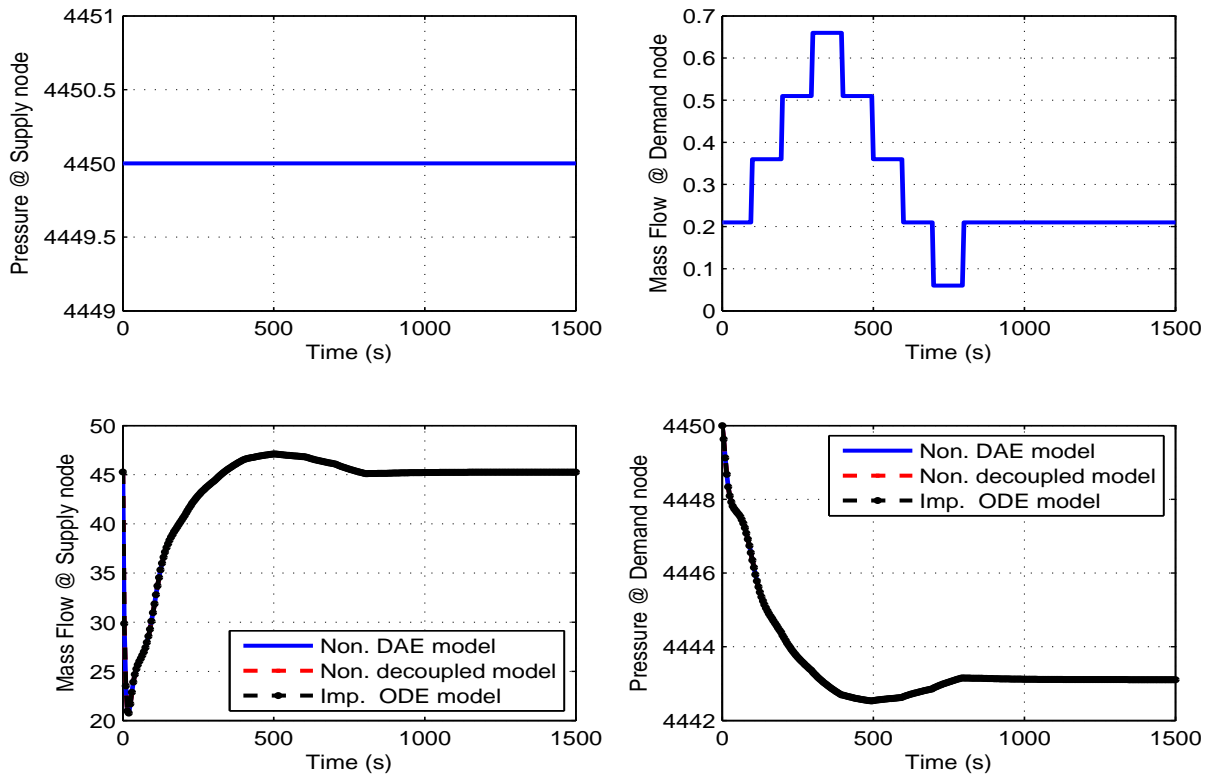}
%
%
\caption{Comparison of the output solutions  in the time interval $t\in \left[0,{1500}s \right].$}
\label{Banagaaya_fig3:1}       
\end{figure}
\end{example}

\begin{example}
\colr{
 In this example, we compare the  matrix properties  of the matrix  pencils  of the derived models and the values of the nonlinear term at a fixed state vector.  In Figures  \ref{Fig_5s:1}- \ref{Fig_5s:3}, we compare the sparsity of the 
 matrix pencils of   the coupled model, decoupled model and implicit ODE model. We can observe that all models are sparse, however the decoupled model is the least sparse.  In Table \ref{Table:5s}, we compare the  finite spectrum of the matrix pencils 
 and the nonlinearity. We can observe that all models have the same spectrum  with   purely imaginary finite  eigenvalues and approximately the same  values of the nonlinear function. 
 \begin{figure}[!h]
 \centering
  \includegraphics[height=.35\textwidth]{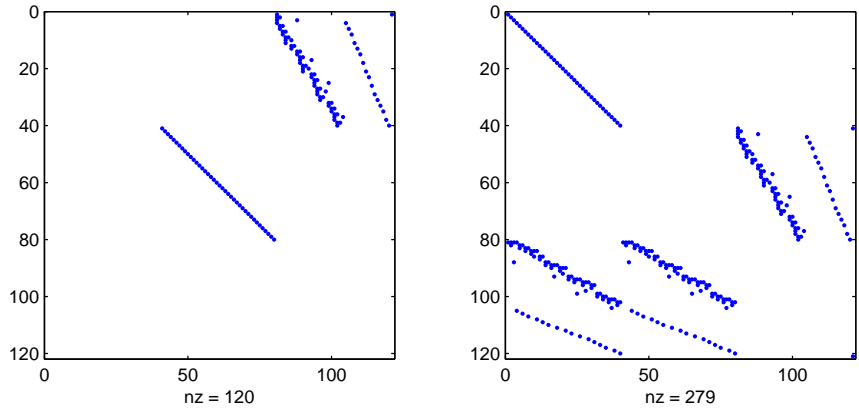} 
  \caption{Sparsity of the matrix pencil $(\E,\A)$ of the coupled model. }
  \label{Fig_5s:1}
 \end{figure}
 \begin{figure}[!h]
 \centering
  \includegraphics[height=.35\textwidth]{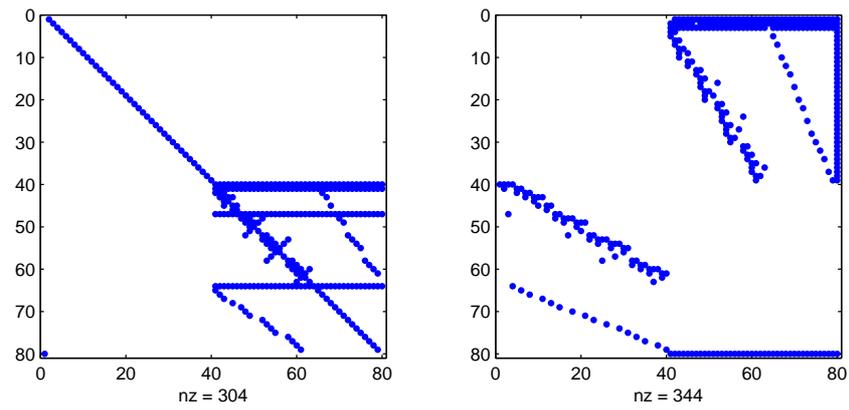} 
 \caption{Sparsity of the matrix pencil $(\E_p,\A_p)$  of the decoupled model. }
  \label{Fig_5s:2}
 \end{figure}
 \begin{figure}[!h]
 \centering
  \includegraphics[height=.35\textwidth]{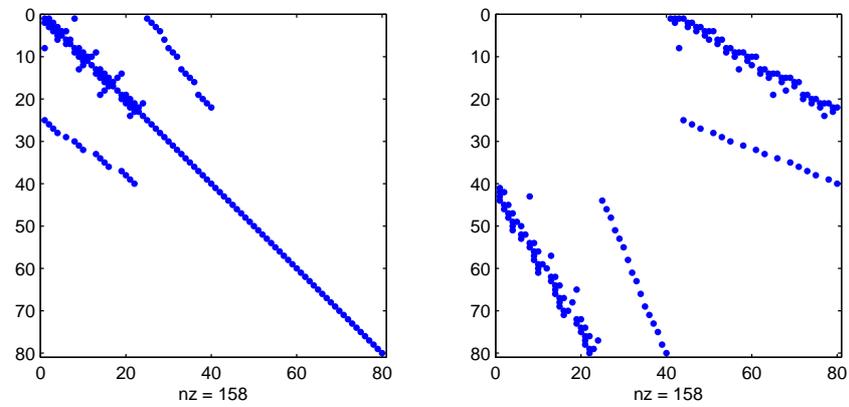} 
 \caption{Sparsity of the matrix pencil of the implicit ODE  model. }
  \label{Fig_5s:3}
 \end{figure}

 \begin{figure}[!h]
 \centering
  \includegraphics[scale=.5]{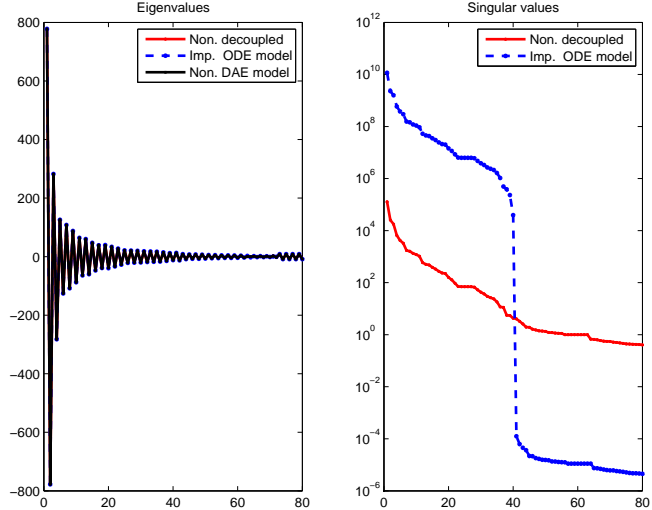} 
  \caption{Comparison of the  eigenvalues and  singular values. }
  \label{Fig_5s:4}
 \end{figure}

 \newpage
 \ \\
  \colr{In Figure \ref{Fig_5s:4},  we compare the  values of the purely imaginary eigenvalues and singular values for  different models. We can observe that eigenvalues exponentially decay for all models. However, the ODE and decoupled models have different 
  singular values. }

 \begin{table}[!h]
   \caption{Comparison of the eigenvalues of the matrix pencil  and the  norm of the nonlinear term}
   \label{Table:5s}
 \begin{center}
 \scalebox{.65}{
\begin{tabular}{|l|l|lll|lll|lll|}
\hline 
$n$ & $n_f$ & \multicolumn{3}{c}{Nonlinear  DAE} &\multicolumn{3}{|c}{ Nonlinear  ODE} &\multicolumn{3}{|c|}{ Nonlinear  Decoupled}\\ \hline 
 & & $\lambda_{min}$ & $\lambda_{max}$ & $\Vert \f(\x)\Vert$ & $\lambda_{min}$ & $\lambda_{max}$ & $\Vert \f(\x)\Vert$ & $\lambda_{min}$ & $\lambda_{max}$ & $\Vert \f(\x)\Vert$\\ \hline 
$4$ & $2$ & $-166.67i$ & $166.67i$ & $2.2202$ &$-166.67i$ & $166.67i$ & $2.2202$ &$-166.67i$ & $166.67i$& $2.2202$\\
$25$ & $16$ & $-0.020803i$ & $1.3352i$ & $0.35903$ &$-0.020803i$ & $1.3352i$  & $0.35903$ & $-0.020803i$ & $1.3352i$  & $0.36265$\\
$55$ & $36$ & $-7.56\times 10^{-4}i$ & $39.558i$ & $69.8168$ &$-7.56\times 10^{-4}i$ & $39.558i$  & $69.8168$ & $-7.56\times 10^{-4}i$ & $39.558i$  & $69.8168$\\
$121$ & $80$ & $-0.4768i$ & $777.7542i$ & $0.53219$ & $-0.4768i$ & $777.7542i$  & $0.53219$ &  $-0.4768i$ & $777.7542i$  & $0.53219$\\ \hline 
 \end{tabular}
 }
 \end{center}
 \end{table}

}
\end{example}

\subsection{Model order reduction}
Here, we illustrate the performance of the proposed IMOR method \colr{compared to  existing MOR methods}.
\begin{example}
\label{Examp:3}
We consider a large-scale   gas transport pipeline network  with $5{,}000$ pipes, $1$ supply node and $1$ demand node. This model was generated numerically using the 
following data.
The length, diameter and average roughness of each pipe are chosen \colp{as} $0.726 \mathrm{m}, \, 1.422\mathrm{m}$ and $1.0\times 10^{-6}\mathrm{m},$ respectively. 
The gas composition is with  specific gas constant $1530 \mathrm{J/KgK}$ at steady pressure  $50\mathrm{bar}$  at supply node and 
mass flow  as a step function as shown in the first row of Figure \ref{Fig_5:1} at the  demand node  at a time interval $t\in \left[0,86400 \right].$
 This lead\colp{s} to a nonlinear  DAE \eqref{Eqn:55} of dimension $n=15,001$. It took $63.7s$ to automatically   decouple the nonlinear DAE in{to} $n_p=10,000$ nonlinear differential equations 
and $n_q=5,001$ algebraic equations. We also generated an index reduced ODE \eqref{ode} of dimension $\tilde{n}=10,000.$  We reduced the decoupled system using POD on both the differential and algebraic parts.

 \begin{table}[!h]
  \caption{Comparison of the ROMs}
  \label{Tab:banagaaya_1}
 \begin{center}
\begin{center}
\begin{tabular}{lrrrr}
\hline 
ROMs & Red. Size ($r$) & \% Red.  & Output error & Speed-ups\\ \hline 
DAE-POD & $2$ & $99.99$  & $3.3\times 10^{-5}$ & $52.9$\\
ODE-POD & $1$ & $9\colp{9}.99$ & $2.1\times 10^{-5}$ & $49.4$\\
I-POD & $6$ & $99.96$  & $1.1\times 10^{-5}$ & $27.0$\\ \hline               \end{tabular}
              \end{center}

 \end{center}
  \label{Tab_5:1}
 \end{table}
 
 \ \\
 Then,  we obtained an I-POD model with  $r_p=2$ and $r_q=4$ leading to a total reduction of $r=r_p+r_q=6\ll 15,001.$
We also used POD to reduce both the nonlinear DAE and ODE directly.  For comparison\colp{,}  the size of ROMs for  different MOR methods is determined by making sure that the output error is below $10^{-4}$  and the results are presented  in  Table \ref{Tab:banagaaya_1}.
All numerical integration was done using implicit-Euler method \colp{with} a fixed time step $h=250$ and LU \colp{based} numerical solver was used for linear solving. 
We can observe that  I-POD leads to the largest ROM \colr{ and lowe\colp{st} speed-ups.} \colp{T}his is due to the fact that its ROM is a DAE  while \colp{the} other ROMs  are  ODEs.
The comparison of the  mass flow at the supply node and  the pressure at the demand node of all ROMs \colp{are}  shown in Figure \ref{Fig_5:1}. 

\begin{figure}[!h]
 \centering
  \includegraphics[width=\textwidth,height=.65\textwidth]{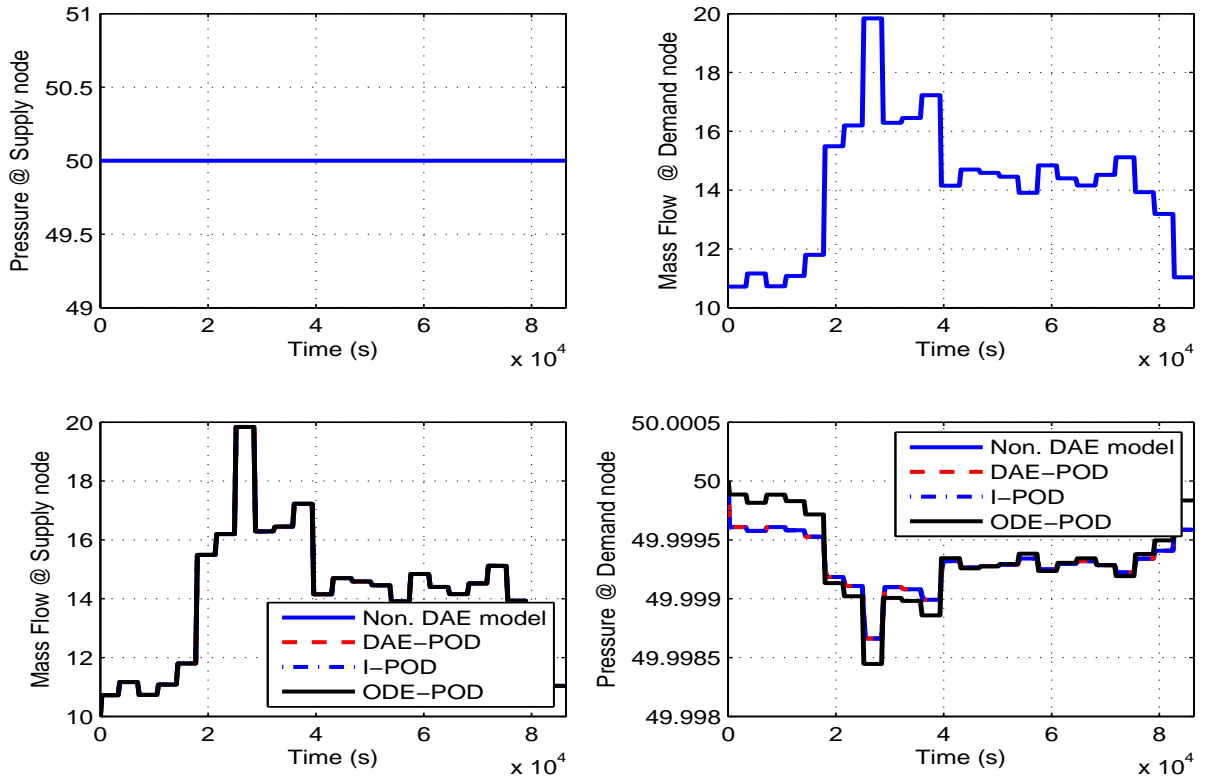} 
  \caption{Comparison of the  pressure at demand nodes and mass flow at supply node. }
  \label{Fig_5:1}
 \end{figure}

 \ \\
 In Figure \ref{Fig_5:2}, we compare \colp{the} output relative error for pressure and mass flow  for different sizes of ROMs. We  can observe \colb{that}  I-POD  is the most accurate    while  ODE-POD is  the le\colp{ast} accurate.
 However, I-POD leads to a s\colp{l}ightly bigger ROM.
 \begin{figure}[!h]
 \centering
  \includegraphics[width=\textwidth,height=.65\textwidth]{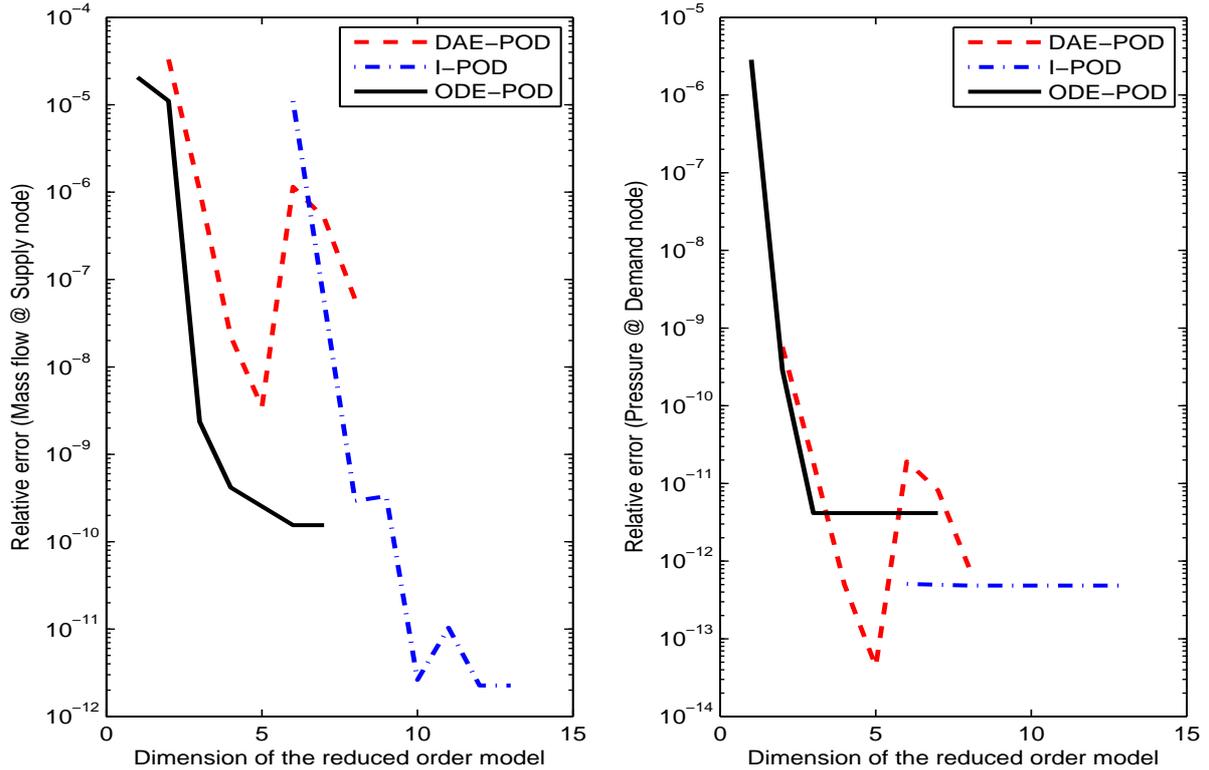} 
  \caption{Comparison of the relative error of the ROMs. }
  \label{Fig_5:2}
 \end{figure}
\end{example}

%
%
\section{Conclusions}
We have proposed a new automatically decoupling strategy and an IMOR method  for nonlinear DAEs with a special nonlinear term. 
This approach eliminates the index problem during simulation and MOR which allows the use of standard numerical integration
methods and MOR techniques. We have derived both the implicit \eqref{eqn:decp111} and explicit \eqref{eqn:303} decoupled systems for index $1$ nonlinear DAEs.  
We have demonstrated  the accuracy  of this approach by applying it \colp{to} 
nonlinear DAEs arising from the gas transportation networks. 
\colp{The computational cost of this approach can be improved by applying  reordering algorithms after decoupling. }
 However, we have restricted ourselves \colb{to} nonlinear DAEs of tractability index one. Future research will deal with  nonlinear DAEs of tractability  index greater than one.

 \bibliographystyle{plain}
\bibliography{reference}

\begin{thebibliography}{10}

\bibitem{Bana2014}
G.~Al\`i, N.~Banagaaya, W.H.A. Schilders, and C.~Tischendorf.
\newblock Index-aware model order reduction for differential-algebraic
  equations.
\newblock {\em Mathematical and Computer Modelling of Dynamical Systems},
  20(4):345--373, 2014.

\bibitem{morAnt05}
A.~C. Antoulas.
\newblock {\em Approximation of Large-Scale Dynamical Systems}.
\newblock {SIAM} Publications, Philadelphia, PA, 2005.

\bibitem{Bana:2014}
N.~Banagaaya.
\newblock {\em {I}ndex-aware model order reduction methods}.
\newblock PhD thesis, Eindhoven University of Technology, Eindhoven,
  Netherlands, 2014.

\bibitem{BanaBook}
N.~Banagaaya, G.~Al\`i, and W.H.A. Schilders.
\newblock {\em {I}ndex-aware {M}odel {O}rder {R}eduction {M}ethods:
  {A}pplications to {D}ifferential-{A}lgebraic {E}quations}, volume~2 of {\em
  Atlantis Studies in Scientific Computing in Electromagnetics}.
\newblock Atlantis Press, 2016.

\bibitem{Bana2018}
N.~Banagaaya, P.~Benner, L.~Feng, P.~Meuris, and W.~Schoenmaker.
\newblock {An index-aware parametric model order reduction method for
  parameterized quadratic differential-algebraic equations}.
\newblock {\em Applied Mathematics and Computation}, 319(C):409--424, 2018.

\bibitem{morBauBF14}
U.~Baur, P.~Benner, and L.~Feng.
\newblock Model order reduction for linear and nonlinear systems: A
  system-theoretic perspective.
\newblock {\em Arch. Comput. Methods Eng.}, 21(4):331--358, 2014.

\bibitem{morBenGW15}
P.~Benner, S.~Gugercin, and K.~Willcox.
\newblock A survey of model reduction methods for parametric systems.
\newblock {\em SIAM Review}, 57(4):483--531, 2015.

\bibitem{MarzBook}
E.~Griepentrog and R.~M\"{a}rz.
\newblock {\em Differential-algebraic equations and their numerical treatment}.
\newblock Teubner, Leipzig, 1986.

\bibitem{morGruHKetal13}
S.~Grundel, N.~Hornung, B.~Klaassen, P.~Benner, and T.~Clees.
\newblock Computing surrogates for gas network simulation using model order
  reduction.
\newblock In S.~Koziel and L.~Leifsson, editors, {\em Surrogate-Based Modeling
  and Optimization}, pages 189--212. Springer, New York, 2013.

\bibitem{morGruHR16}
S.~Grundel, N.~Hornung, and S.~Roggendorf.
\newblock Numerical aspects of model order reduction for gas transportation
  networks.
\newblock In S.~Koziel, L.~Leifsson, and X.-S. Yang, editors, {\em
  Simulation-Driven Modeling and Optimization}, pages 1--28. Springer, 2016.

\bibitem{morGruJ15}
S.~Grundel and L.~Jansen.
\newblock Efficient simulation of transient gas networks using {IMEX}
  integration schemes and {MOR} methods.
\newblock In {\em 54th IEEE Conference on Decision and Control (CDC), Osaka,
  Japan}, pages 4579--4584, December 2015.

\bibitem{morGruJHetal14}
S.~Grundel, L.~Jansen, N.~Hornung, T.~Clees, C.~Tischendorf, and P.~Benner.
\newblock Model order reduction of differential algebraic equations arising
  from the simulation of gas transport networks.
\newblock In {\em Progress in Differential-Algebraic Equations},
  Differential-Algebraic Equations Forum, pages 183--205. Springer Berlin
  Heidelberg, 2014.

\bibitem{Canonical}
G.~Kalogeropoulos, M.~Mitrouli, A.~Pantelous, and D.~Triantafyllou.
\newblock The {W}eierstra{\ss} canonical form of a regular matrix pencil:
  Numerical issues and computational techniques.
\newblock In S.~Margenov, L.G. Vulkov, and J.~Wa\'{s}niewski, editors, {\em
  Numerical Analysis and Its Applications}, pages 322--329. Springer, Berlin,
  Heidelberg, 2009.

\bibitem{MehrBook}
P.~Kunkel and V.~Mehrmann.
\newblock {\em Differential Algebraic Equations: Analysis and Numerical
  Solution}, volume~1.
\newblock EMS, 2006.

\bibitem{Marz92}
R.~M\"{a}rz.
\newblock Numerical methods for differential algebraic equations.
\newblock {\em Acta Numerica}, 21(5):141--198, 1992.

\bibitem{Marz4}
R.~M{\"a}rz.
\newblock Canonical projectors for linear differential algebaric equations.
\newblock {\em Computers Math. Applications}, 31(4/5):121--135, 1996.

\bibitem{Marz02}
R.~M{\"a}rz.
\newblock The index of linear differential algebraic equations with properly
  stated leading terms.
\newblock {\em Results in Math.}, 42:308--338, 2002.

\bibitem{Marzlead}
R.~M\"{a}rz.
\newblock Solvability of linear differential algebraic equations with properly
  stated leading terms.
\newblock {\em Results in Math.}, 45(1):88--105, 2004.

\bibitem{Marz_arrays2}
R.~M\"{a}rz.
\newblock Characterizing differential algebraic equations without the use of
  derivative arrays.
\newblock {\em J. Comput. Math. Appl.}, 50(7):1141--1156, 2005.

\bibitem{4lectures}
Steffen Schulz.
\newblock Four {L}ectures on {D}ifferential-{A}lgebraic {E}quations.
\newblock Research Report 497, The University of Auckland, Department of
  Mathematics, June 2003.

\bibitem{Kaarthik18}
K.~Sundar and A.~Zlotnik.
\newblock State and parameter estimation for natural gas pipeline networks
  using transient state data.
\newblock {\em {IEEE} Trans. Control Syst.}, pages 1--15, 2018.

\bibitem{Zhang}
Z.~Zhang and N.~Wong.
\newblock An efficient projector-based passivity test for descriptor systems.
\newblock {\em IEEE Trans. Comput.-Aided Design Integr. Circuits Syst.},
  29(8):1203--1214, Aug 2010.

\end{thebibliography}

\end{document}